\documentclass{amsproc}

\usepackage{verbatim} 
\usepackage{amsmath}
\usepackage{amsfonts}
\usepackage{amssymb}
\usepackage{mathrsfs}
\usepackage{amsthm}
\usepackage{newlfont}

\usepackage{fullpage,amsmath, amssymb,amscd}
\usepackage{latexsym, epsfig, color}
\usepackage[mathscr]{eucal}
\usepackage{xypic, graphicx}
\usepackage{amscd}
\usepackage[all, cmtip]{xy}
\usepackage{tikz,tikz-cd,xcolor}
\usepackage{tikz-cd}

\usepackage{float, hyperref, nccmath, multirow, caption, subcaption}

\newcommand\cyr{%
 \renewcommand\rmdefault{wncyr}%
 \renewcommand\sfdefault{wncyss}%
 \renewcommand\encodingdefault{OT2}%
\normalfont\selectfont} \DeclareTextFontCommand{\textcyr}{\cyr}

\newtheorem{theorem}{Theorem}
\newtheorem{lemma}[theorem]{Lemma}
\newtheorem{corollary}[theorem]{Corollary}
\newtheorem{proposition}[theorem]{Proposition}
\newtheorem{remark}[theorem]{Remark}

\def\Z{\mathbb Z}

\def\Q{\mathbb Q}
\def\R{\mathbb R}
\def\C{\mathbb C}
\def\F{\mathbb F}

\def\G{\mathbb G}

\def\Z{\mathbb Z}

\def\Q{\mathbb Q}
\def\R{\mathbb R}
\def\C{\mathbb C}
\def\F{\mathbb F}

\def\G{\mathbb G}

\def\Im{\operatorname{Im}}
\def\Re{\operatorname{Re}}

\def\det{\operatorname{det}}

\def\tr{\operatorname{tr}}

\def\End{\operatorname{End}}

\def\Gal{\operatorname{Gal}}

\def\Frob{\operatorname{Frob}}

\def\id{\operatorname{id}}
\def\mod{\operatorname{mod}}

\def\disc{\operatorname{disc}}
\def\car{\operatorname{char}}

\def\GL{\operatorname{GL}}

\def\SL{\operatorname{SL}}

\def\li{\operatorname{li}}

\def\O{\operatorname{O}}
\def\o{\operatorname{o}}

\def\log{\operatorname{log}}

\def\ds{\displaystyle}

\begin{document}

\title{
Bounds for the distribution of the Frobenius traces associated to products of non-CM elliptic curves
}

\date{}

\author{
Alina Carmen Cojocaru and Tian Wang}
\address[Alina Carmen  Cojocaru]{
\begin{itemize}
\item[-]
Department of Mathematics, Statistics and Computer Science, University of Illinois at Chicago, 851 S Morgan St, 322
SEO, Chicago, 60607, IL, USA;
\item[-]
Institute of Mathematics  ``Simion Stoilow'' of the Romanian Academy, 21 Calea Grivitei St, Bucharest, 010702,
Sector 1, Romania
\end{itemize}
} \email[Alina Carmen  Cojocaru]{cojocaru@uic.edu}

\address[Tian Wang]{
\begin{itemize}
\item[-]
Department of Mathematics, Statistics and Computer Science, University of Illinois at Chicago, 851 S Morgan St, 322
SEO, Chicago, 60607, IL, USA;
\end{itemize}
} \email[Tian Wang]{twang213@uic.edu}

\renewcommand{\thefootnote}{\fnsymbol{footnote}}
\renewcommand{\thefootnote}{\arabic{footnote}}

\renewcommand{\thefootnote}{\fnsymbol{footnote}}
\renewcommand{\thefootnote}{\arabic{footnote}}

\thanks{A.C.C. was partially supported  by a Collaboration Grant for Mathematicians from the Simons Foundation  
under Award No. 709008. }

\thanks{February 9, 2022}

\begin{abstract}
Let $g \geq 1$ be an integer
and 
let $A/\Q$ be an abelian variety that is
isogenous over $\Q$ to 
a product of $g$  elliptic curves defined over $\Q$, 
pairwise  non-isogenous over $\overline{\Q}$
and
each without complex multiplication.
For an  integer $t$ and a positive real number $x$, 
denote by $\pi_A(x, t)$ the number of primes $p \leq x$, 
of good reduction for 
$A$,
 for which the Frobenius trace $a_{1, p}(A)$ associated to the reduction of $A$ modulo $p$ equals $t$.
Assuming the Generalized Riemann Hypothesis for Dedekind zeta functions, 
we prove that
$\pi_A(x, 0) \ll_A x^{1 - \frac{1}{3 g+1 }}/(\log x)^{1 - \frac{2}{3 g+1}}$
and 
$\pi_A(x, t) \ll_A x^{1 - \frac{1}{3 g + 2}}/(\log x)^{1 - \frac{2}{3 g + 2}}$ if $t \neq 0$.
These bounds largely improve upon recent ones obtained for $g = 2$ by H. Chen, N. Jones, and V. Serban,
 and may be viewed as generalizations to arbitrary $g$ of the bounds obtained for $g=1$  by  M.R. Murty, V.K. Murty, and N. Saradha,  combined with a refinement in the power of $\log x$ by D. Zywina.
Under the same assumptions, we also prove the existence of a density one set of primes $p$  satisfying 
$|a_{1, p}(A)|>p^{\frac{1}{3 g + 1} - \varepsilon}$ for any  fixed $\varepsilon>0$.
\end{abstract}

\maketitle


\section{Introduction}

A prominent open problem in arithmetic geometry, formulated by S. Lang and H. Trotter in the 1970s 
(\cite[Part I]{LaTr76}),
concerns the distribution of the Frobenius traces associated to the reductions modulo primes of an  elliptic curve defined over 
$\Q$ and without complex multiplication. 
In recent years,  this problem has been formulated in broader settings, such as that of abelian varieties 
(e.g., \cite{CoDaSiSt17}, \cite{ChJoSe20},  and \cite{Ka09}).  
The goal of the present article is to provide upper bounds related to the distribution of the Frobenius traces defined by the product of  non-isogenous
elliptic curves  defined over $\Q$ and having no complex multiplication, as explained below.

Let $g \geq 1$ be an integer
and 
let $A/\Q$ be an abelian variety that is
isogenous over $\Q$ to 
a product of $g$  elliptic curves defined over $\Q$, 
pairwise 
non-isogenous over $\overline{\Q}$
and
each without complex multiplication.
Denote by $N_A$ the conductor of $A$.
For a prime $p \nmid N_A$, 
we write the characteristic polynomial of the Frobenius endomorphism 
acting on the reduction of $A$ modulo $p$ as
$$
P_{A, p}(X) =
 X^{2g} + a_{1, p}(A) X^{2g-1} + \ldots  + a_{2g-1, p}(A) X + a_{2g, p}(A)  \in \Z[X].
$$
Then, for a fixed integer $t$, we study the counting function
$$
\pi_A(x, t) := \#\left\{p \leq x: p \nmid N_A,  a_{1, p}(A) = t\right\}.
$$

When $g = 1$, 
the asymptotic behavior of $\pi_A(x, t)$ is predicted by the Lang-Trotter Conjecture on Frobenius traces (see \cite[Part I, p. 33]{LaTr76} for the original and \cite[Conjectures 2-3, p. 435]{BaJo09} for refined versions).
Specifically, in \cite{LaTr76}, it is conjectured that 
there  exists an explicit non-negative constant $C(A, t)$, 
which depends on the arithmetic of $A$ and  $t$, such that, as $x \rightarrow \infty$,
$$
\pi_A(x, t) \sim C(A, t) \frac{\sqrt{x}}{\log x}.
$$
When $g = 2$ and $t \neq 0$, a similar asymptotic formula is predicted by a conjecture of  H. Chen, N. Jones and V. Serban (see \cite[Conjecture 1.0.5, p. 3]{ChJoSe20}). 
An analogous asymptotic behavior is expected to hold for arbitrary $t$ and $g$ in our given setting.

When $g =1$, M.R. Murty, V.K. Murty and N. Saradha \cite[Theorem, p. 254]{MuMuSa88} built on the work of 
J-P. Serre \cite[Section 8, pp. 188-191]{Se81}
and
proved that, under the assumption of a Generalized Riemann Hypothesis,
$\pi_{A}(x, 0) \ll_{A} x^{\frac{3}{4}}$
and 
$\ds \pi_{A}(x, t) \ll_{A} \frac{x^{\frac{4}{5}}}{(\log x)^{\frac{1}{5}}}$ if $t \neq 0$.
More recently, D. Zywina \cite[Theorem 1.2, p. 236]{Zy15} obtained 
the improvements 
$\ds \pi_A(x, 0) \ll_{A} \frac{x^{\frac{3}{4}}}{(\log x)^{\frac{1}{2}}}$
and
$\ds \pi_A(x, t) \ll_{A} \frac{x^{\frac{4}{5}}}{(\log x)^{\frac{3}{5}}}$ if $t \neq 0$.
The bound 
$\pi_A(x, 0) \ll_{A} x^{\frac{3}{4}}$ 
was also proven
unconditionally
using
 the interpretation of the condition
$a_{1, p}(A) = 0$ as $p$ being a supersingular prime for the elliptic curve $A$ 
(see \cite[Theorem B, p. 131]{El91}).
When $t \neq 0$, only the  bound $\ds \pi_A(x, t) \ll_{A} \frac{x (\log \log x)^2}{(\log x)^2}$ is known unconditionally, 
thanks to the work of J. Thorner and A. Zaman \cite[Theorem 1.4, p. 4997]{ThZa17},
 who built on prior results of 
J-P. Serre \cite[Section 8, pp. 188-191]{Se81},
D. Wan \cite[Theorem 1.3, p. 250]{Wan90}, 
and 
V.K. Murty \cite[Theorem 5.1, p. 302]{Mu97}.

When $g = 2$ and $t \neq 0$,
H. Chen, N. Jones and V. Serban \cite[Theorem 2.4.1, p. 14]{ChJoSe20} proved  that, 
under the assumption of a Generalized Riemann Hypothesis, we have
$\ds \pi_A(x, t) \ll_{A} \frac{x^{\frac{13}{14}}}{(\log x)^{\frac{5}{7}}}$;
unconditionally,
they proved that
$\ds \pi_A(x, t) \ll_{A} \frac{x (\log \log x)^{\frac{2}{7}} (\log \log \log x)^{\frac{1}{7}} }{(\log x)^{\frac{8}{7}}}$.

\medskip

In this paper, we prove conditional upper bounds 
for $\pi_A(x, t)$
that largely improve upon the  conditional bound of H. Chen, N. Jones, and V. Serban when $g = 2$, 
and recover the conditional  bound of M.R. Murty, V.K. Murty and N. Saradha and D. Zywina when $g=1$.

\begin{theorem}\label{main-thm1}
Let $t \in \Z$
and
let $A/\Q$ be an abelian variety that is isogenous over $\Q$ to 
a product of $g$  elliptic curves defined over $\Q$, 
pairwise non-isogenous over $\overline{\Q}$
and
each without complex multiplication.
Under the assumption 
of Generalized Riemann Hypothesis  (GRH) for Dedekind zeta functions,
we have that, for  any sufficiently large $x$,
$$
\pi_A(x, t) 
\ll_{A}
\left\{
\begin{array}{cl}
 \ds\frac{x^{1 - \frac{1}{3 g+1}}}{(\log x)^{1 - \frac{2}{3 g+1}}}   & \text{if $t = 0$,}
\\
\\
\ds\frac{x^{1 - \frac{1}{3 g + 2}}}{ (\log x)^{1 - \frac{2}{3 g + 2}} } & \text{if $t \neq 0$.}
\end{array}
\right.
$$
\end{theorem}

\medskip

As an immediate application of Theorem  \ref{main-thm1},
we obtain a result about the non-lacunarity of the sequence $(a_{1, p}(A))_{p}$,
that is, about the non-zero values of $a_{1, p}(A)$. 
 In fact, we obtain a result about any fixed value of $a_{1, p}(A)$.
\begin{corollary}\label{main-corollary1}
In the setting and under the assumptions of Theorem \ref{main-thm1},  
we have that, 
as $x \rightarrow \infty$,
\begin{equation}\label{corollary-1}
\#\left\{p \leq x: p \nmid N_A, a_{1, p}(A) \neq t \right\} 
\sim \pi(x).
\end{equation}
\end{corollary}

As another quick application of Theorem \ref{main-thm1},  we obtain
 the existence of a density one  set of primes $p$ with large Frobenius traces $a_{1, p}(A)$.
Namely, in the setting and under the assumptions of Theorem \ref{main-thm1}, we have that, 
for any $\varepsilon > 0$ and
as $x \rightarrow \infty$,
\begin{equation}\label{corollary-2}
\hspace*{-0.5cm}
\#\left\{p \leq x: 
p \nmid N_A,
 |a_{1, p}(A)| > p^{\frac{1}{3 g + 2} - \varepsilon} \right\} 
\sim \pi(x).
\end{equation}

Without applying Theorem \ref{main-thm1} directly,  but instead
using a proof strategy similar to that of the proof of Theorem \ref{main-thm1}, 
we can improve upon (\ref{corollary-2}) to show the following.
\begin{theorem}\label{main-thm2}
Let $A/\Q$ be an abelian variety that is isogenous over $\Q$ to 
a product of $g$  elliptic curves defined over $\Q$, 
pairwise  non-isogenous over $\overline{\Q}$
and
each without complex multiplication.
Under the assumption 
of Generalized Riemann Hypothesis  (GRH) for Dedekind zeta functions,
we have that,
for any $\varepsilon > 0$ and
as $x \rightarrow \infty$,
$$
\#\left\{p \leq x: 
p \nmid N_A,
 |a_{1, p}(A)| > p^{\frac{1}{3 g + 1} - \varepsilon} \right\} 
\sim \pi(x).
$$
\end{theorem}

The above results may be re-written solely in terms of  a $g$-tuple of
elliptic curves $E_1/\Q, \ldots, E_g/\Q$,
assumed to be pairwise non-isogenous over $\overline{\Q}$
and each  without complex multiplication,
 as follows.
Denote by $N_{E_1}, \ldots, N_{E_g}$ the conductors of $E_1, \ldots, E_g$, respectively.
For each  integer $1 \leq i \leq g$ and for each prime $p \nmid N_{E_i}$, 
 denote by $a_p(E_i)$ the integer  defined by
$|\overline{E}_i(\F_p)| =p + 1 - a_p(E_i)$,
where 
 $\overline{E}_i$  is the
 reduction of $E_i$ modulo $p$.
Taking $A := E_1 \times \ldots \times E_g$,
we have  that $a_{1, p}(A) = - (a_p(E_1) + \ldots + a_p(E_g))$
for any prime $p \nmid N_{E_1} \ldots N_{E_g}$ (to be explained in Section 6).
Now let $t \in \Z$.
Then equivalent formulations of
Theorem \ref{main-thm1}, Corollary \ref{main-corollary1}, and Theorem \ref{main-thm2}
are that, 
under the assumption
of GRH,
for any sufficiently large $x$, we have
$$
\#\left\{
p \leq x:
p \nmid N_{E_1} \ldots N_{E_g},
a_p(E_1) + \ldots + a_p(E_g) =  t
\right\}
\ll_{E_1, \ldots, E_g}
\left\{
\begin{array}{cl}
 \ds\frac{x^{1 - \frac{1}{3 g+1}}}{(\log x)^{1 - \frac{2}{3 g+1}}}   & \text{ {\emph{if}} $t = 0$,}
\\
\\
\ds\frac{x^{1 - \frac{1}{3 g + 2}}}{ (\log x)^{1 - \frac{2}{3 g + 2}} } & \text{ {\emph{if}} $t \neq 0$,}
\end{array}
\right.
$$

\medskip
and 
for any $\varepsilon > 0$ and
as $x \rightarrow \infty$,
we have

\begin{equation*}
\#\left\{
p \leq x:
p \nmid N_{E_1} \ldots N_{E_g},
a_p(E_1) + \ldots + a_p(E_g) \neq t
\right\}
\sim \pi(x),
\end{equation*}

\begin{equation*}
\#\left\{
p \leq x:
p \nmid N_{E_1} \ldots N_{E_g},
\left|a_p(E_1) + \ldots + a_p(E_g)\right| > p^{\frac{1}{3 g + 1} - \varepsilon}
\right\}
\sim \pi(x).
\end{equation*}

\bigskip

The general strategy of proving 
an upper bound for $\pi_A(x, t)$
 is to relax the equality $a_{1, p}(A) = t$ to
a congruence $\tr \bar{\rho}_{A, m} \left(\left(\frac{\Q(A[m])/\Q}{p}\right)\right)  \equiv - t (\mod m)$,
where
 $m$ is an arbitrary integer coprime to $p$,
$\bar{\rho}_{A, m}$ is the Galois representation associated to 
the $m$-division field  $\Q(A[m])$ of $A$,
and
$\left(\frac{\Q(A[m])/\Q}{p}\right)$
is the Artin symbol at $p$ in the Galois group $\Gal(\Q(A[m])/\Q)$.
An upper bound for $\pi_A(x, t)$
may then be derived by interpreting this congruence as a Chebotarev condition
in the finite Galois extension $\Q(A[m])/\Q$ and by invoking an effective version of the Chebotarev Density Theorem. 

Inspired by \cite{MuMuSa88}, 
we  pursue a refinement of this general strategy 
and relax the equality $a_{1, p}(A) = t$ to a  Chebotarev condition that holds in a suitably chosen proper subextension of $\Q(A[m])/\Q$ in which Artin's Holomorphy Conjecture is known to hold. 
When $g \geq 2$, we unravel a few such possible subextensions and  carry out the refined strategy in the subextension that leads to the best  result, in particular the subextension that leads to a generalization 
of the best upper bounds known for $g=1$.

When studying the Frobenius traces of an elliptic curve, it is natural to place a special focus on
elliptic curves $A/\Q$ having the property that
$\End_{\overline{\Q}}(A) \simeq \Z$, since this case is regarded as  generic.
When studying the Frobenius traces of an abelian variety $A$ of higher dimension,
it is then natural to place an initial focus on the following two cases:  
that of an abelian variety $A/\Q$ with $\End_{\overline{\Q}}(A)  \simeq \Z$, as pursued in \cite{CoDaSiSt17} and in an upcoming paper 
by the present authors,
and
that of an abelian variety
that is $\Q$-isogenous 
to a product of $g$ elliptic curves defined over $\Q$,
pairwise  non-isogenous over $\overline{\Q}$
and
each without complex multiplication,
as pursued in the present paper.
Note that, when $A/\Q$ is  $\Q$-isogenous 
to a product of $g$ elliptic curves defined over $\Q$
that are pairwise isogenous over $\overline{\Q}$, the study of $\pi_A(x, t)$ for $g \geq 2$ reduces to the study of 
$\pi_A(x, t)$ for $g = 1$. Other complementary cases  remain to be pursued separately. 

Finally, let us note that other variations of questions inspired by the Lang-Trotter Conjecture have been investigated in the setting of pairs of elliptic curves in works such as
\cite{FoMu95},  \cite{AkDaJu04}, and \cite{AkPa19}.

\bigskip

{\bf{General notation.}}
Throughout the paper, we use the following  notation. 

\medskip

$\bullet$
For a set $S$, we denote its cardinality by $|S|$ or $\# S$.

\medskip

$\bullet$
Given suitably defined real functions $h_1, h_2$,
we say that 
$h_1 = \o(h_2)$ if $\ds\lim_{x \rightarrow \infty} \frac{h_1(x)}{h_2(x)} = 0$;
we say that
$h_1 = \O(h_2)$ or, equivalently, that $h_1 \ll h_2$, 
if 
$h_2$ is positive valued 
and
 there exists a positive constant $c$ such that 
$|h_1(x)| \leq c \ h_2(x)$ for all $x$ in the domain of $h_1$ and $h_2$;
we say that
$h_1 \asymp h_2$ 
if 
$h_1$, $h_2$ are positive valued 
and 
$h_1 \ll h_2 \ll h_1$;
we say that
$h_1 = \O_D(h_2)$ or, equivalently, that $h_1 \ll_D h_2$,
if 
$h_1 = \O(h_2)$
and
the implied $\O$-constant $c$ depends on priorly given data  $D$;
we say that 
$h_1 \asymp_D h_2$ 
if the implied constant $c$ in 
at least
 one of the $\ll$-bounds
$h_1 \ll h_2 \ll h_1$
depends on priorly  given data $D$;
finally, we say that 
$h_1 \sim h_2$ 
if
$\ds\lim_{x \rightarrow \infty} \frac{h_1(x)}{h_2(x)} = 1$.

\medskip

$\bullet$
We use the letters $p$ and $\ell$ to denote rational primes.
We use $\pi(x)$ to denote the number of primes $p \leq x$,
 and
$\li x$ to denote the logarithmic integral $\ds\int_{2}^x \frac{1}{\log t} \ d t$. 
Note that $\ds \li x \sim \frac{x}{\log x}$
and 
recall that the Prime Number Theorem asserts that
$\pi(x) \sim \li x$.
 
 \medskip

$\bullet$
 Given a number field $K$, 
 we denote  by  ${\cal{O}}_K$ its ring of integers;
 we denote by $\overline{K}$ a fixed algebraic closure of $K$;
 we denote by $\Gal(\overline{K}/K)$ the absolute Galois group of $K$;
 we denote by ${\sum}_K$ the set of non-zero prime ideals of $K$;
 we denote by $[K:\Q]$ the degree of $K$ over $\Q$;
 we denote by $d_K \in \Z$ the discriminant of an integral basis of ${\cal{O}}_K$
 and by  $\disc(K/\Q) = \Z d_K \unlhd \Z$ the discriminant ideal of $K/\Q$.
  For a prime ideal $\mathfrak{p} \in {\sum}_K$, 
 we denote by $N_{K/\Q}(\mathfrak{p})$ its norm in $K/\Q$
 and by $\Frob_{\mathfrak{p}} \in \Gal(\overline{K}/K)$ its Frobenius class.
 We say that $K$ satisfies the Generalized Riemann Hypothesis (GRH) if
 the Dedekind zeta function $\zeta_K$ of $K$ has the property that,
 for any $\rho \in \C$ with $0 \leq \Re \rho \leq 1$ and $\zeta_K(\rho) = 0$, we have $\Re(\rho) = \frac{1}{2}$. When $K=\Q$, the Dedekind zeta function is the Riemann zeta function, in which case GRH is typically referred to as RH.
 
 \medskip
 
 $\bullet$
 For a nonzero unitary commutative ring $R$, we denote by
 $R^{\times}$  its group of multiplicative units.
For an integer $n \geq 1$, 
we denote by
$I_n$ the  identity $n \times n$ matrix with entries in $R$.
For an arbitrary $n \times n$ matrix $M$ with entries in $R$,
we denote by $\tr M$ and $\det M$ its trace and determinant.
We define the general linear group  $\GL_{n}(R)$
 as the collection of $n \times n$ matrices $M$ with entries in $R$ and with $\det M \in R^{\times}$.
 
 \medskip
 
 $\bullet$
 For a positive integer $m$, we denote by $\Z/m \Z$ the ring of integers modulo $m$.
 For an arbitrary prime $\ell$, we denote by $\Z_{\ell}$ the ring of $\ell$-adic integers.
We set $\hat{\Z} := \ds \lim_{\leftarrow \atop m} \Z/m \Z$
and note that there is a ring isomorphism 
$\hat{\Z}  \simeq \ds\prod_{\ell} \Z_{\ell}$.

\medskip

\noindent
{\bf{Acknowledgments.}} We thank Jacob Mayle for a careful reading of our manuscript and for useful remarks. We are deeply grateful to the referees for their detailed comments and suggestions.


\bigskip
\section{The Chebotarev Density Theorem}


The proofs of  Theorems \ref{main-thm1} and \ref{main-thm2} are based on multi-step
 applications of different effective versions of the Chebotarev Density Theorem, which we now recall.

Let $L/K$ be a  Galois extension of number fields, with $G:= \Gal(L/K)$, and let
${\cal{C}} \subseteq G$ be a union of conjugacy classes of $G$; throughout the paper, we will assume ${\cal{C}}\neq \emptyset$. 
 We denote by $[L:K]$ the degree of $L$ over $K$ and  by $\disc(L/K) \unlhd {\cal{O}}_K$ the discriminant ideal of $L/K$.
 For each $\mathfrak{p} \in {\sum}_K$,
 we denote 
 by $\wp \in {\sum}_L$ an arbitrary prime ideal with $\wp \mid \mathfrak{p}$,
 by ${\cal{D}}_{\wp}$ the decomposition group of $\wp$ in $L/K$,
 and
 by  ${\cal{I}}_{\wp}$ the inertia group of $\wp$ in $L/K$.
 We denote
 by $ \left(\frac{L/K}{\wp}\right) \in {\cal{D}}_{\wp}/{\cal{I}}_{\wp}$ the Artin symbol associated to $\wp$,
and 
by 
$\left(\frac{L/K}{\mathfrak{p}}\right)  
:=  
\left\{\gamma\in G: 
\exists \wp\mid \mathfrak{p} \text{ such that } \gamma\in {\cal{D}}_{\wp}  \ \text{and} \ \gamma {\cal{I}}_{\wp} =\left(\frac{L/K}{\wp}\right) 
\right\}$
the  Artin symbol associated to $\mathfrak{p}$
in the extension $L/K$.
 Notice that, for any integer $m \geq 1$,
 the set
$\left(\frac{L/K}{\mathfrak{p}}\right)^m
:=
\left\{
\gamma^m: 
\gamma \in \left(\frac{L/K}{\mathfrak{p}}\right)\right\}
 \subseteq G$
 is  invariant under conjugation by elements of $G$.
 We set
 $$
 {\cal{P}}(L/K) := \{p: \: \ \exists \ \mathfrak{p} \in {\sum}_{K}
 \ \text{such that} \
 \mathfrak{p} \mid p \ \text{and} \ \mathfrak{p} \mid \disc(L/K)\},
 $$
 $$
 M(L/K) := 2 [L:K] |d_K|^{\frac{1}{[K:\Q]}} \ds\prod_{p \in {\cal{P}}(L/K)} p,
 $$
 and recall that
\begin{equation}\label{hensel}
\log \left|N_{K/\Q} (\disc (L/K))\right|
\leq
([L : \Q] - [K : \Q])
\left(
\ds\sum_{p \in {\cal{P}}(L/K)} \log p
\right)
+
[L : \Q] \log [L : K]
\end{equation}
(see \cite[Proposition 5, p.~129]{Se81}).
 
 The extension $L/K$ is said to satisfy the Artin Holomorphy Conjecture (AHC) if,
 for any irreducible character $\chi: G \longrightarrow \C$,
 the Artin $L$-function $L(s, \chi)$ extends to a function that is analytic on $\C$ when $\chi \neq 1$
 and that is analytic on $\C \backslash \{1\}$ when $\chi = 1$. It is known that AHC holds if $\Gal(L/K)$ is abelian (\cite{Ar27}).

We denote by $\delta_{\cal{C}} : G \longrightarrow \{0, 1\}$ the characteristic function on ${\cal{C}}$.
For  any  prime ideal $\mathfrak{p} \in {\sum}_K$,
 we choose an arbitrary prime ideal $\wp \in {\sum}_L$ with $\wp \mid \mathfrak{p}$,
and then, 
for any integer $m \geq 1$,
we define
$$
\delta_{\cal{C}}\left(\left(\frac{L/K}{\mathfrak{p}}\right)^m\right)
:=
\frac{1}{|{\cal{I}}_{\wp}|} 
\ds\sum_{
\gamma \in {\cal{D}}_{\wp}
\atop{
\gamma {\cal{I}}_{\wp} = \left(\frac{L/K}{\wp}\right)^m 
}
}
\delta_{\cal{C}}(\gamma).
$$
Note that this definition is independent of the choice of $\wp$.
For an arbitrary real number $x > 0$, we define
\begin{eqnarray*}
\pi_{\cal{C}}(x, L/K)
&:=&
\ds\sum_{
\mathfrak{p} \in {\sum}_K
\atop{
\mathfrak{p} \nmid \disc(L/K)
\atop{
N_{K/\Q}(\mathfrak{p})\leq x
}
}
}
\delta_{\cal{C}}\left(
\left(\frac{L/K}{\mathfrak{p}}\right)
\right),
\\
\widetilde{\pi}_{\cal{C}}(x, L/K)
&:=&
\ds\sum_{m \geq 1}
\frac{1}{m}
\ds\sum_{
\mathfrak{p} \in {\sum}_K
\atop{
N_{K/\Q}(\mathfrak{p}^m) \leq x
}
}
\delta_{\cal{C}}\left(
\left(\frac{L/K}{\mathfrak{p}}\right)^m
\right).
\end{eqnarray*}
For simplicity of notation, 
we use
$\pi_1(x, L/K)$
and
$\widetilde{\pi}_1(x, L/K)$
to denote 
$\pi_{\{\id_L\}}(x, L/K)$
and
$\widetilde{\pi}_{\{\id_L\}}(x, L/K)$, respectively. 

The Chebotarev Density Theorem states that, as $x \rightarrow \infty$,
\begin{equation}\label{chebotarev}
\pi_{\cal{C}}(x, L/K) \sim \frac{|{\cal{C}}|}{|G|} \li x.
\end{equation}
Since the growth in $x$ of the counting function $\widetilde{\pi}_{\cal{C}}(x, L/K)$ is closely related to that of $\pi_{\cal{C}}(x, L/K)$
thanks to the bound
\begin{equation}\label{pi-versus-pi-tilde}
\left|
\widetilde{\pi}_{\cal{C}}(x, L/K) - \pi_{\cal{C}}(x, L/K)
\right|
\ll
[K:\Q]
\left(
\frac{x^{\frac{1}{2}}}{\log x} + \log M(L/K)
\right)
\end{equation}
(see \cite[Proposition 7, p. 138]{Se81}), 
asymptotic (\ref{chebotarev}) is equivalent to 
\begin{equation}\label{chebotarev-tilde}
\widetilde{\pi}_{\cal{C}}(x, L/K) \sim \frac{|{\cal{C}}|}{|G|} \li x.
\end{equation}
The advantage of considering $\widetilde{\pi}_{\cal{C}}(x, L/K)$ over $\pi_{\cal{C}}(x, L/K)$ is that the former satisfies the following two functorial properties.

\begin{proposition}\label{functorial-prop}
Let $L/K$ be a Galois extension of number fields with $G:= \Gal(L/K)$ and let
${\cal{C}} \subseteq G$ be a union of conjugacy classes of $G$.

\begin{enumerate}
\item[(i)]\label{functorial-subgroup}
Let $H \leq G$ be a subgroup of $G$ such that
every element of ${\cal{C}}$ is conjugate to some element of $H$.
Denote by $L^H$ the subfield of $L$ fixed by $H$.
Then, for any $x > 0$,
$$
\widetilde{\pi}_{\cal{C}}(x, L/K)
\leq
\widetilde{\pi}_{{\cal{C}} \cap H}(x, L/L^H).
$$
\item[(ii)]\label{functorial-quotient-group}
Let $N \unlhd G$ be a normal subgroup of $G$ such that $N {\cal{C}} \subseteq {\cal{C}}$.
Denote by $L^N$ the subfield of $L$ fixed by $N$.
Denote by $\widehat{\cal{C}} $ the image  of ${\cal{C}}$ in the quotient group $G/N$, viewed as $\Gal(L^N/K)$.
Then, for any $x > 0$,
$$
\widetilde{\pi}_{\cal{C}}(x, L/K) = \widetilde{\pi}_{\widehat{\cal{C}}}(x, L^N/K).
$$
\end{enumerate}
\end{proposition}
\begin{proof}
See \cite[Proposition 8, pp. 139--140]{Se81}.
\end{proof}

\begin{corollary}\label{functorial-upper-bound}
Let $L/K$ be a  Galois extension of number fields with $G:= \Gal(L/K)$, and let
${\cal{C}} \subseteq G$ be a union of conjugacy classes of $G$.
Let $H \leq G$ be a subgroup of $G$ such that
every element of ${\cal{C}}$ is conjugate to some element of $H$.
Let $N \unlhd H$ be a normal subgroup of $H$
such that
$N ({\cal{C}} \cap H)\subseteq {\cal{C}} \cap H$. 
Denote by $\widehat{\cal{C} \cap H} $ the image  of ${\cal{C}} \cap H$ in the quotient group $H/N$, viewed as $\Gal(L^N/L^H)$.
Then, for any $x > 0$,
$$
\widetilde{\pi}_{\cal{C}}(x, L/K) \leq \widetilde{\pi}_{\widehat{\cal{C} \cap H}}\left(x, L^N/L^H\right).
$$
Consequently, for any sufficiently large $x$,
\begin{eqnarray}\label{pi-versus-pi-tilde-upper}
\pi_{\cal{C}}(x, L/K)
&\ll&
\pi_{\widehat{\cal{C} \cap H}}\left(x, L^N/L^H\right)
\nonumber
\\
&&
+
[K:\Q]
\left(
\frac{x^{\frac{1}{2}}}{\log x} + \log M(L/K)
\right)
+
[L^H:\Q]
\left(
\frac{x^{\frac{1}{2}}}{\log x} + \log M(L^N/L^H)
\right).
\end{eqnarray}
\end{corollary}
\begin{proof}
The corollary  follows by combining parts (i) and (ii) of Proposition \ref{functorial-prop} with   (\ref{pi-versus-pi-tilde}).
\end{proof}

In spite of its strength, the asymptotic formula  (\ref{chebotarev}) for $\pi_{{\cal{C}}}(x, L/K)$ does not suffice for our purposes.
Instead,  the proofs of our main results will require an effective version  with an explicit error term, such as the following conditional result.

\begin{theorem}\label{cheb-LaOd} 
Let $L/K$ be a Galois extension of number fields, with $G:= \Gal(L/K)$, and let
${\cal{C}} \subseteq G$ be a union of conjugacy classes of $G$.
Assume that GRH holds for the number field $L$.
 Then there exists an absolute  constant
 $c > 0$ such that, for any $x>2$,
 \begin{eqnarray*}
\left|
\pi_{\cal{C}}(x, L/K) - \frac{|\cal{C}|}{|G|} \pi(x)
\right|
\leq
c \frac{|\cal{C}|}{|G|} x^{\frac{1}{2}}
\left(
\log |d_L|
+
[L : \Q] \log x
\right).
\end{eqnarray*}
 \end{theorem}
 \begin{proof}
The original reference is  \cite{LaOd77}. For this variation,  see \cite[Th\'eor\`eme~4, p.~133]{Se81}.
 \end{proof}

The proofs of our main results will also use the following 
conditional  upper bound for $\pi_{\cal{C}}(x, L/K)$.
\begin{theorem}\label{chebotarev-upper-bound}
Let $L/K$ be a  Galois extension of number fields with $G:= \Gal(L/K)$, and let
${\cal{C}} \subseteq G$ be a union of conjugacy classes of $G$.
Assume that GRH holds for the number field $L$ and that AHC holds for the number field extension $L/K$.
Then,  for any sufficiently large $x$,
$$
\pi_{\cal{C}}(x, L/K)
\ll
\frac{|\cal{C}|}{|G|}\cdot  \frac{x}{\log x}
+
|\cal{C}|^{\frac{1}{2}} [K:\Q] \left(\log M(L/K)\right) \cdot \frac{x^{\frac{1}{2}}}{\log x}.
$$
\end{theorem}
\begin{proof}
See
\cite[Theorem 2.3, p. 240]{Zy15}. 
\end{proof}
\begin{corollary}\label{cheb-upper-bound-as-wanted}
Let $L/K$ be a  Galois extension of number fields with $G:= \Gal(L/K)$, and let
${\cal{C}} \subseteq G$ be a union of conjugacy classes of $G$.
Let $H \leq G$ be a subgroup of $G$ such that 
every element of ${\cal{C}}$ is conjugate to some element of $H$.
Let $N \unlhd H$ be a normal subgroup of $H$
such that
$N ({\cal{C}} \cap H)\subseteq {\cal{C}} \cap H$. 
Denote by $\widehat{\cal{C} \cap H} $ the image  of ${\cal{C}} \cap H$ in the quotient group $H/N$, viewed as $\Gal(L^N/L^H)$.
Assume that the group $H/N$ is abelian
and that  GRH holds for $L^N$. Then, for any sufficiently large $x$,
\begin{eqnarray*}
\pi_{\cal{C}}(x, L/K)
&\ll&
\frac{\left|\widehat{\cal{C} \cap H}\right| \cdot |N|}{|H|} \cdot \frac{x}{\log x}
\\
&+&
\left|\widehat{{\cal{C}} \cap H}\right|^{\frac{1}{2}} [L^H :\Q] \frac{x^{\frac{1}{2}}}{\log x} \log M(L^N/L^H)
\\
&+&
[K:\Q] \left(\frac{x^{\frac{1}{2}}}{\log x} + \log M(L/K)\right)
+
[L^H:\Q] \left(\frac{x^{\frac{1}{2}}}{\log x} + \log M(L^N/L^H)\right).
\end{eqnarray*}
\end{corollary}
\begin{proof}
First use  (\ref{pi-versus-pi-tilde-upper}), and then Theorem \ref{chebotarev-upper-bound}  applied to the extension $L^N/L^H$.
\end{proof}

\bigskip
\section{An application of the effective version of the Chebotarev Density Theorem}

The proofs of Theorems \ref{main-thm1} and \ref{main-thm2} will make use of the following application of 
Theorem \ref{cheb-LaOd}.
\begin{lemma}\label{max-lemma}
Let ${\cal{S}}$ be a non-empty set of rational primes,
let $(K_p)_{p \in {\cal{S}}}$ be a family of finite Galois extensions of $\Q$,
and
let $(\cal{C}_p)_{p\in \cal{S}}$ be a family of non-empty sets such that
 each $\cal{C}_p$  is a union  of conjugacy classes of $\Gal(K_p/\Q)$.
Assume that
there exist an absolute constant $c_1 > 0$
and 
a  function $f: \R \to (0, \infty)$
such that

\begin{equation}\label{hypothesis1}
\ds
[K_p:\Q]
\leq
c_1
\quad \text{for any sufficiently large} \ p,
\end{equation}


\begin{equation}\label{hypothesis0}
\ds \log |d_{K_p}| 
\leq
 f(z)
\quad \text{for any sufficiently large } z  \text{ and all 
} p\leq z.
\end{equation}
 For each 
 $x>2$, let $y = y(x) > 2$, $u =u(x) > 2$ be such that
\begin{equation}\label{hypothesis3}
u \leq y,
\end{equation}
and assume that, for any $\varepsilon > 0$,
\begin{equation}\label{hypothesis4}
u \geq 
c_2(\varepsilon) y^{\frac{1}{2}} (\log y)^{2 + \varepsilon} 
\
\text{for some constant} \ 
c_2(\varepsilon) > 0
\end{equation}
and
\begin{equation}\label{hypothesis5}
\ds\lim_{x \rightarrow \infty} \frac{f(x)}{(\log y)^{1 + \varepsilon}} = 0.
\end{equation}
Assume GRH for Dedekind zeta functions. 
Then, for any $\varepsilon > 0$,
there exists a constant $c(\varepsilon) > 0$ such that, for any sufficiently large $x$,
\begin{equation}\label{max-lemma-bound}
\#\left\{p \leq x: p \in {\cal{S}}\right\}
\leq
c(\varepsilon)
\ds\max_{
y \leq \ell \leq y + u
} 
\#\left\{p \leq x: p \in {\cal{S}}, \ell \nmid d_{K_p}, \left(\frac{K_p/\Q}{\ell} \right)\subseteq \cal{C}_p \right\}.
\end{equation}
\end{lemma}

\begin{proof}
The proof follows the strategy of  \cite[Lemma 4.4, p. 269]{MuMuSa88}, adapted to our general setting.
We start by fixing $x > 2$, $y = y(x) > 2$, $u = u(x) > 2$, and $\varepsilon > 0$ 
such that (\ref{hypothesis3}) - (\ref{hypothesis5}) hold.
Then, we observe that
\begin{eqnarray}
&&
(\pi(y + u) - \pi(y))
\ds\max_{y \leq \ell \leq y + u} 
\#\left\{p \leq x: p \in {\cal{S}}, \ell \nmid d_{K_p}, \left(\frac{K_p/\Q}{\ell} \right)\subseteq \cal{C}_p \right\}
\label{max-lower}
\\
&\geq&
\ds\sum_{y \leq \ell \leq y  + u}
\#\left\{p \leq x: p \in {\cal{S}}, \ell \nmid d_{K_p}, \left(\frac{K_p/\Q}{\ell} \right)\subseteq \cal{C}_p \right\}
\nonumber
\\
&=&
\ds\sum_{
p \leq x
\atop{
p \in {\cal{S}}
}
}
\#\left\{y \leq \ell \leq y + u: \ell \nmid d_{K_p}, \left(\frac{K_p/\Q}{\ell} \right)\subseteq \cal{C}_p 
\right\}
\nonumber
\\
&=&
\ds\sum_{
p \leq x
\atop{
p \in {\cal{S}}
}
}
\left(\pi_{\cal{C}_p}(y+u, K_p/\Q) - \pi_{\cal{C}_p}(y, K_p/\Q)\right).
\nonumber
\end{eqnarray}

To estimate the above difference, we apply Theorem \ref{cheb-LaOd} under GRH and deduce that
\begin{eqnarray}\label{difference-pi1}
&&
\ds\sum_{
p \leq x
\atop{
p \in {\cal{S}}
}
}
\left(\pi_{\cal{C}_p}(y+u, K_p/\Q) - \pi_{\cal{C}_p}(y, K_p/\Q)\right)
\\
&=&
\ds\sum_{
p \leq x
\atop{p \in {\cal{S}}}
}
\frac{|\cal{C}_p|}{[K_p:\Q]}
(\pi(y + u) - \pi(y))
+
\ds\sum_{
p \leq x
\atop{p \in {\cal{S}}}
}
\left(
E_1(y, u, \cal{C}_p, K_p)
+
E_2(y, u, \cal{C}_p, K_p)
\right)
\nonumber
\end{eqnarray}
for some real-valued functions 
$E_1(y, u, \cal{C}_p, K_p)$, $E_2(y, u,\cal{C}_p, K_p)$, which depend on $y$, $u$, $\cal{C}_p$, and $K_p$,
and for which there exist absolute constants 
$c_3 > 0$ and $c_4 > 0$ such that
\begin{equation*}\label{bound-E1-p}
|E_1(y, u,\cal{C}_p, K_p)|
\leq
c_3 (y + u)^{\frac{1}{2}}
\left(
\frac{|\cal{C}_p| \log |d_{K_p}|}{[K_p : \Q]} + |\cal{C}_p| \log (y + u)
\right),
\end{equation*}

\begin{equation*}\label{bound-E2-p}
|E_2(y, u, \cal{C}_p, K_p)|
\leq
c_4 y^{\frac{1}{2}}
\left(
\frac{|\cal{C}_p| \log |d_{K_p}|}{[K_p : \Q]} + |\cal{C}_p|  \log y
\right).
\end{equation*}
Recalling assumptions (\ref{hypothesis1}) - 
 (\ref{hypothesis0}), we obtain the upper bound
\begin{eqnarray}\label{bound-E1-plus-E2}
\left|
\ds\sum_{
p \leq x
\atop{p \in {\cal{S}}}
}
\left(
E_1(y, u, \cal{C}_p, K_p) 
+
E_2(y, u, \cal{C}_p, K_p)
\right)
\right|
&\leq&
c_5  (y + u)^{\frac{1}{2}}(f(x) + \log(y+u)) \ \#\left\{p \leq x: p \in {\cal{S}}\right\}
\end{eqnarray}
for some  absolute constant $c_5 > 0$.

Thanks to
(\ref{hypothesis4})
and the assumption of GRH (specifically, RH in this case), there exists  some constant 
$c_6(\varepsilon) > 0$  such that
\begin{equation}\label{RH-difference-pi}
\pi(y + u) - \pi(y) \geq 
c_6(\varepsilon)  \frac{u}{\log u}.
\end{equation}
In particular, $\pi(y + u) - \pi(y) > 0$. 

Putting together (\ref{max-lower}) - (\ref{RH-difference-pi})
and
dividing the resulting inequality by $\pi(y + u) - \pi(y)$, 
we obtain that
\begin{eqnarray}\label{max-greater-than}
&&
\ds\max_{y \leq \ell \leq y + u} 
\#\left\{p \leq x: p \in {\cal{S}}, \ell \nmid d_{K_p}, \left(\frac{K_p/\Q}{\ell} \right)\subseteq \cal{C}_p\right\}
\nonumber
\\
&\geq&
\ds\sum_{
p \leq x
\atop{p \in {\cal{S}}}
}
\frac{|\cal{C}_p|}{[K_p:\Q]}
+
\O_{\varepsilon}
\left(
\frac{\log u}{u}
(y + u)^{\frac{1}{2}}
(f(x) + \log (y + u))
\
\#\{p \leq x: p \in {\cal{S}}\}
\right).
\end{eqnarray}
Invoking assumptions 
(\ref{hypothesis3}) - (\ref{hypothesis4}),
we deduce that the above $\O_{\varepsilon}$-term becomes
$$
\O_{\varepsilon}
\left(
\frac{f(x) + \log y}{(\log y)^{1 + \varepsilon}}
\
\#\{p \leq x: p \in {\cal{S}}\}
\right).
$$
Invoking assumption (\ref{hypothesis5}), we deduce that the  latter $\O_{\varepsilon}$-term is
$\o(\#\{p \leq x: p \in {\cal{S}}\})$.
Finally, recalling (\ref{hypothesis1}) and that each $\cal{C}_p$ is non-empty, 
we deduce that
\begin{equation*}
\ds\sum_{
p \leq x
\atop{p \in {\cal{S}}}
}
\frac{|\cal{C}_p|}{[K_p:\Q]}
\geq
\frac{1}{c_1}
\
\#\{p \leq x: p \in {\cal{S}}\}.
\end{equation*}
Then
(\ref{max-greater-than}) gives
\begin{equation*}
\ds\max_{y \leq \ell \leq y + u} 
\#\left\{p \leq x: p \in {\cal{S}}, \ell \nmid d_{K_p}, \left(\frac{K_p/\Q}{\ell} \right)\subseteq \cal{C}_p\right\}
\gg_{\varepsilon}
\#\{p \leq x: p \in {\cal{S}}\},
\end{equation*}
as desired.
\end{proof}

\begin{remark}\label{max-lemma-hypotheses}
{\emph{
In the proof of Lemma \ref{max-lemma},
instead of hypotheses (\ref{hypothesis1}) - 
(\ref{hypothesis0}), it suffices to make the 
weaker hypotheses that
there exist absolute constants $c'_1, c'_2, c'_3 > 0$ 
such that,
for any sufficiently large $z$, 
\begin{equation*}\label{hypothesis1bis}
\ds\sum_{p \leq z  \atop{ p \in {\cal{S}}}}
\frac{\left| \cal{C}_p \right|}{[K_p : \Q]}
\geq
c'_1 \ \#\{p \leq z: p \in {\cal{S}}\},
\end{equation*}
\begin{equation*}\label{hypothesis2bis}
\ds\sum_{p \leq z  \atop{ p \in {\cal{S}}}}
\frac{\left| \cal{C}_p \right| \log |d_{K_p}|}{[K_p : \Q]}
\leq
c'_2 f(z) \ \#\{p \leq z: p \in {\cal{S}}\},
\end{equation*}
\begin{equation*}\label{hypothesis0bis}
\ds\sum_{p \leq z  \atop{ p \in {\cal{S}}}}
\left|\cal{C}_p\right| 
\leq
c'_3 \ \#\{p \leq z: p \in {\cal{S}}\}.
\end{equation*}.
}}
\end{remark}

\bigskip
%
\section{Subgroups of  $\GL_{2}(\Z/\ell \Z)^g$}

Our main contribution to generalizing the method of M.R. Murty,  V.K. Murty, and N. Saradha \cite{MuMuSa88} to the case of a product
of $g$ non-isogenous elliptic curves without complex multiplication
 consists of
unravelling suitable Galois extensions $L/K$ of number fields with $\Gal(L/K) < \GL_{2}(\Z/\ell\Z)^g$ for some well-chosen prime $\ell$,
and 
of unravelling suitable conjugacy classes $\cal{C} \subseteq \Gal(L/K)$
for which we can invoke the version of the Chebotarev Density Theorem given in
Corollary \ref{cheb-upper-bound-as-wanted}.
With this goal in mind, we devote Section 4 to investigations of particular subgroups of $\GL_{2}(\Z/\ell \Z)^g$.

For a fixed rational prime $\ell$, our main focus will be on the subgroup
\begin{equation}\label{group-G-of-ell}
G(\ell)
:= \{(M_1, \ldots, M_g)\in \G(\ell): \det M_1 = \ldots = \det M_g\}
\end{equation}
of the product group
$$
\G(\ell) :=\GL_2(\ell)^g,
$$
where 
$$\GL_2(\ell) := \GL_2(\Z/\ell\Z).$$
Denoting
by ${\cal{G}}_m(\ell)$ the multiplicative group $(\Z/\ell \Z)^{\times}$
and
by ${\cal{G}}_a(\ell)$ the additive group $\Z/\ell \Z$,
we set

\medskip
$
 \mathcal{B}_{\GL_2}(\ell)
  :=
   \left\{
  \begin{pmatrix} 
  a_1  & b_1\\ 0 & a_2  
  \end{pmatrix}
  \in \GL_2(\ell) : 
  a_1, a_2 \in \mathcal{G}_m(\ell),
  b_1\in \mathcal{G}_a(\ell) \right\},
 $
 
 \medskip
  $
 \mathcal{U}_{\GL_2}(\ell)
  :=
   \left\{
  \begin{pmatrix} 
  1 & b\\ 0 & 1  
  \end{pmatrix}
  \in \GL_2(\ell) : 
  b\in \mathcal{G}_a(\ell) \right\},
 $

 \medskip
  $
 \mathcal{T}_{\GL_2}(\ell)
  :=
   \left\{
  \begin{pmatrix} 
  a_1 & 0\\ 0 & a_2  
  \end{pmatrix}
  \in \GL_2(\ell) : 
 a_1, a_2 \in \mathcal{G}_m(\ell) \right\},
 $

\medskip
$
B(\ell) := \{(M_1, \ldots, M_g)\in  \mathcal{B}_{\GL_2}(\ell)^g: \det M_1 = \ldots = \det M_g\},
$

\medskip
$
U(\ell) := \mathcal{U}_{\GL_2}(\ell)^g,
$

\medskip
$
U'(\ell) := (\Z/\ell\Z)^{\times}\cdot U(\ell),
$

\medskip
$
T(\ell) :=  \{(M_1, \ldots, M_g)\in  \mathcal{T}_{\GL_2}(\ell)^g: \det M_1 = \ldots = \det M_g\}.
$



\noindent

\medskip
The following lemmas will be crucial in the proof of our main theorems.
 
  \begin{lemma}\label{centr-G-B-U}  
  For any prime $\ell$,
  we have that
   $U(\ell)$ and  $U'(\ell)$ are normal subgroups of $B(\ell)$, 
   and that $B(\ell)/U(\ell)$ are $B(\ell)/U'(\ell)$ are abelian groups.
 \end{lemma}

 \begin{proof}
The subgroup $\mathcal{U}_{\GL_2}(\ell)$ is normal in $\mathcal{B}_{\GL_2}(\ell)$, 
and
the quotient  $\mathcal{B}_{\GL_2}(\ell)/\mathcal{U}_{\GL_2}(\ell)$  is abelian, since
it satisfies a group isomorphism
 $\mathcal{B}_{\GL_2}(\ell)/\mathcal{U}_{\GL_2}(\ell) \simeq  \mathcal{T}_{\GL_2}(\ell)$.
 Consequently, the subgroup $U(\ell)$ is normal in $\mathcal{B}_{\GL_2}(\ell)^g$, with 
   $
   \mathcal{B}_{\GL_2}(\ell)^g/U(\ell) 
   $
   an abelian group. By the inclusion of subgroups
   $
   U(\ell) \subseteq U'(\ell) \subseteq B(\ell) \subseteq \mathcal{B}_{\GL_2}(\ell)^g,
   $
   we deduce
    that 
  both $U(\ell)$ and $U'(\ell)$ are normal subgroups in $B(\ell)$ 
   and that
 both  $B(\ell)/U(\ell)$ and $B(\ell)/U'(\ell)$ are abelian.
 \end{proof}


\begin{lemma}\label{counting-G(l)-etc}
For any prime $\ell$, we have
\begin{eqnarray*}
|\GL_2(\ell)| 
&=&
(\ell - 1) \ell (\ell^2 - 1),
\\
|\G(\ell)|
&=&
(\ell - 1)^g \ell^g (\ell^2 - 1)^g,
\\
|\mathcal{B}_{\GL_2}(\ell)^g|
&=&
(\ell - 1)^{2 g} \ell^{ g},
\\
|G(\ell)| 
&=& 
(\ell - 1) \ell^{g} (\ell^2 - 1)^{g},
 \\
 |B(\ell)| 
 &=&
 (\ell-1)^{g + 1} \ell^g,
 \\
 |U(\ell)|
 &=&
 \ell^{g},
 \\
 \left|U'(\ell)\right|
 &=&
 (\ell - 1) \ell^g,
 \\
 |T(\ell)|
 &=& 
 (\ell-1)^{g + 1}.
\end{eqnarray*}
\end{lemma}
\begin{proof}

The formula for the order of $\GL_2(\ell)$ is well-known and easy to derive.
The  formulae for the orders of  
$\G(\ell)$, $\mathcal{B}_{\GL_2}(\ell)^g$, $U(\ell)$, $U'(\ell)$, and $T(\ell)$ are clear from the definitions of the groups. 
To find the orders of $G(\ell)$ and $B(\ell)$, observe that there are short exact sequences of finite groups
$$
1 \to G(\ell) \to \G(\ell) \xrightarrow{\det_{1, \ldots, g-1}}  \mathcal{G}_m(\ell)^{g-1} \to 1,
$$
$$
1 \to B(\ell) \to \mathcal{B}_{\GL_2}(\ell)^g \xrightarrow{\det_{1, \ldots, g-1}}   \mathcal{G}_m(\ell)^{g-1} \to 1,
$$
where
$\det_{1, \ldots, g-1}(M_1, \ldots, M_g) := ((\det M_1) (\det M_2)^{-1}, \ldots, (\det M_1) (\det M_g)^{-1})$.
We deduce that
$$|G(\ell)| = |\G(\ell)|/|\mathcal{G}_m(\ell)|^{g-1} = (\ell-1) \ell^g (\ell^2 -1)^g,$$
$$|B(\ell)| = |\mathcal{B}_{\GL_2}(\ell)^g|/|\mathcal{G}_m(\ell)|^{g-1}  = (\ell-1)^{g+1} \ell^g.$$
\end{proof}

\begin{lemma}\label{B/U-lemma}
For any prime $\ell$,
we have a group isomorphism $B(\ell)/U(\ell) \simeq T(\ell)$.
\end{lemma}
\begin{proof}
Consider the  group homomorphisms
$
T(\ell) \hookrightarrow B(\ell) \twoheadrightarrow B(\ell)/U(\ell).
$
The composition has kernel $T(\ell)\cap U(\ell)= 1$, hence it is injective.
Since, by Lemma \ref{counting-G(l)-etc},  each of the groups  $B(\ell)/U(\ell)$ and $T(\ell)$ has size $(\ell-1)^{g+1}$,
the composition must be an isomorphism.
\end{proof}

\bigskip
\section{Conjugacy classes of subgroups of $G(\ell)$}

Let $\ell \nmid 2 g$ be a fixed prime and let $t \in \Z$ be a fixed integer.
We devote Section 5 to investigations of particular 
unions of conjugacy classes of
$G(\ell)$, $B(\ell)$, and $T(\ell)$.
As usual, when we regard $\Z/\ell \Z$ as a field, we use the notation $\F_{\ell}$. 
In this case, we denote by 
$\overline{\F}_{\ell}$ a fixed algebraic closure of $\F_{\ell}$.

Before introducing the unions of conjugacy classes we are interested in investigating, 
let us observe that
 if $M = (M_1, \ldots, M_g) \in G(\ell)$,
then the characteristic polynomial of $M$
relates to
 the characteristic polynomials of 
 $M_1, \ldots, M_g \in \GL_2(\ell)$ 
 through the equation
 \begin{equation}\label{char-poly-M-M1-Mg}
\car_M(X)=\car_{M_1}(X)\cdot \ldots \cdot \car_{M_g}(X)  \in  \F_{\ell}[X]. 
\end{equation}
When working over the algebraic closure $\overline{\F}_{\ell}$ of the field $\F_{\ell}$, 
for each $1 \leq i \leq g$
we may write the linear factorization of $\car_{M_i}(X)$ as
 \begin{equation}\label{char-poly-M1-Mg-factor}
\car_{M_i}(X)
= 
(X - \lambda_1(M_i)) (X - \lambda_2(M_i))
 \in \overline{\F}_{\ell}[X].
\end{equation}
Then, by putting together (\ref{char-poly-M-M1-Mg}) and (\ref{char-poly-M1-Mg-factor}), we see that
 the linear factorization of $\car_{M}(X)$ over $\overline{\F}_{\ell}$ is
\begin{equation}\label{char-poly-M-M1-Mg-factor}
\car_{M}(X)
= 
\ds\prod_{1 \leq i \leq g} (X - \lambda_1(M_i)) (X - \lambda_2(M_i))
 \in \overline{\F}_{\ell}[X].
\end{equation}

\medskip

Using the above notation and fixing $z > 0$, 
we now introduce the sets to be investigated in this section: 
\begin{eqnarray*}
\cal{C}(\ell, t)
&:=&
\left\{
(M_1, \ldots, M_g) \in G(\ell):
\lambda_1(M_i), \lambda_2(M_i) \in \F_{\ell}^{\times} \ \forall 1 \leq i \leq g,
\ds\sum_{
1 \leq i \leq g
}
\tr M_i
\equiv
- t (\mod \ell)
\right\},
\end{eqnarray*}
\begin{eqnarray*}
\cal{C}_{\text{Borel}}(\ell, t) 
&:=& 
{\cal{C}}(\ell, t) \cap B(\ell),
\\
\widehat{\cal{C}}_{\text{Borel}}(\ell, t) 
&:=&
\; \text{the image of $\cal{C}_{\text{Borel}}(\ell, t)$ in  $B(\ell)/U(\ell)$},
\\
\widehat{\cal{C}}'_{\text{Borel}}(\ell, t)
&:=&
\; \text{the image of $\cal{C}_{\text{Borel}}(\ell, t)$ in  $B(\ell)/U'(\ell)$},
\\
{\cal{C}}_{\text{Torus}}(\ell, t)
&:=&
{\cal{C}}(\ell, t) \cap T(\ell),
\\
%
\cal{C}(\ell, |t|\leq z)  
&:=&
\bigcup_{\substack{t\in \Z\\|t|\leq z}}\cal{C}(\ell, t),
\\
\cal{C}_{\text{Borel}}(\ell, |t|\leq z) 
& := &
\cal{C}(\ell, |t|\leq z) \cap B(\ell),
\\
\widehat{\cal{C}}_{\text{Borel}}(\ell, |t|\leq z) 
& := &
\text{the image of } \cal{C}_{\text{Borel}}(\ell, |t| \leq z) \text{ in } B(\ell)/U(\ell).
\end{eqnarray*}

\medskip

\begin{lemma}\label{propr-conj-sets-first}

\noindent
\begin{enumerate}
\item[(i)]
$\emptyset \neq \cal{C}_{\text{Torus}}(\ell, t) \subseteq  {\cal{C}}_{\text{Borel}}(\ell, t) \subseteq {\cal{C}}(\ell, t)$
and $\cal{C}_{\text{Borel}}(\ell, |t|\leq z)\neq \emptyset$.
\item[(ii)]
$\cal{C}(\ell, t)$ is a union of conjugacy classes of $G(\ell)$.
\item[(iii)]
$\cal{C}_{\text{Borel}}(\ell, t)$ is a union of conjugacy classes of $B(\ell)$.
\item[(iv)]
$\cal{C}_{\text{Torus}}(\ell, t)$ is a  union of conjugacy classes of $T(\ell)$.
\item[(v)]
$U(\ell) \ {\cal{C}}_{\text{Borel}}(\ell, t) \subseteq {\cal{C}}_{\text{Borel}} (\ell, t)$
 and 
$U(\ell) \ \cal{C}_{\text{Borel}}(\ell, |t|\leq z)\subseteq \cal{C}_{\text{Borel}}(\ell, |t|\leq z)$.
\item[(vi)]
$U'(\ell) \ \cal{C}_{\text{Borel}}(\ell, 0) \subseteq \cal{C}_{\text{Borel}} (\ell, 0)$. 
\end{enumerate}
\end{lemma}
\begin{proof}
(i) Since $\ell$ is a prime such $\ell \nmid 2 g$, we know that $(2 g)^{-1} (\mod \ell)$ exists. Now we consider each of the cases $t \neq 0$ and $t = 0$.
When $t \neq 0$,
for each $1\leq i\leq g$,
we take
$$
M_i :=
\begin{pmatrix}
- t (2 g)^{-1} (\mod \ell) &   0 
\\
 0 & - t (2 g)^{-1} (\mod \ell) 
 \\
\end{pmatrix}.
$$
Then 
\[
\lambda_{1}(M_i) \equiv \lambda_2(M_i)\equiv - t (2 g)^{-1} (\mod \ell)  \in \F_{\ell}^{\times} \ \forall 1\leq i\leq g, 
\quad
\sum_{
1 \leq i \leq g
}
\tr M_i
\equiv
- t (\mod \ell),
\]
and so $(M_1, \ldots, M_g) \in
\cal{C}_{\text{Torus}}(\ell, t)$.
When $t = 0$,  
for each $1\leq i\leq g$,
 we take 
$$
M_i :=
\begin{pmatrix}
(2 g)^{-1} (\mod \ell)   & 0
\\
 0  & -(2 g)^{-1} (\mod \ell) 
\end{pmatrix}.
$$
Then 
\[
\lambda_{1}(M_i)\equiv (2 g)^{-1} (\mod \ell)
\
\text{and}
\
 \lambda_2(M_i)\equiv-(2 g)^{-1}(\mod \ell) \in \F_{\ell}^{\times} \ \forall 1\leq i\leq g, \quad 
\sum_{
1 \leq i \leq g
}
\tr M_i
\equiv
0 (\mod \ell),
\]
and so $(M_1, \ldots, M_g) \in
\cal{C}_{\text{Torus}}(\ell, 0)$.
We conclude that $\cal{C}_{\text{Torus}}(\ell, t) \neq \emptyset$.

 The inclusions
$\cal{C}_{\text{Torus}}(\ell, t) \subseteq  {\cal{C}}_{\text{Borel}}(\ell, t) \subseteq {\cal{C}}(\ell, t)$
 are clear from the definitions of the sets.
 
 Finally, observing that 
 $\ds\cal{C}_{\text{Borel}}(\ell, |t|\leq z) 
 = 
 \bigcup_{\substack{t\in \Z\\|t|\leq z}}\cal{C}_{\text{Borel}}(\ell, t)$, 
 the claim that
 $\cal{C}_{\text{Borel}}(\ell, |t|\leq z)\neq \emptyset$  follows 
 from the previous inclusions and the non-emptiness of $\cal{C}_{\text{Borel}}(\ell, t)$.

\medskip
(ii) $\cal{C}(\ell, t)$ is a subset  of $G(\ell)$ characterized by some condition on eigenvalues of $g$ invertible $2\times 2$ matrices. 
Since eigenvalues are invariant under conjugation, $\cal{C}(\ell, t)$ is a  union of conjugacy classes of $G(\ell).$

\medskip
(iii) By (ii), $\cal{C}(\ell, t)$ is invariant under conjugation by $B(\ell)$. Since $B(\ell)$ is obviously invariant under conjugation by $B(\ell)$,  so is the intersection
${\cal{C}}(\ell, t) \cap B(\ell)$. 
This implies that $\cal{C}_{\text{Borel}}(\ell, t)$  is a union of conjugacy classes of $B(\ell)$.
 
\medskip
(iv) Since ${\cal{C}}(\ell, t) $ and $T(\ell)$ are invariant under conjugation by $T(\ell)$,  
the set $\cal{C}_{\text{Torus}}(\ell, t)$ is a union of conjugacy classes of $T(\ell)$.

\medskip
(v) Let 
$M' \in U(\ell)$,
$M \in {\cal{C}}_{\text{Borel}}(\ell, t)$.
Then $M' M \in B(\ell)$, with
$M' M$  having the same diagonal entries as $M$. Since $\tr M = - t$,  
we obtain that $\tr (M' M) = - t$.  
Thus $U(\ell) \ {\cal{C}}_{\text{Borel}}(\ell, t) \subseteq {\cal{C}}_{\text{Borel}} (\ell, t)$.
Upon recalling  that $\ds\cal{C}_{\text{Borel}}(\ell, |t|\leq z) =\bigcup_{\substack{t\in \Z\\|t|\leq z}}\cal{C}_{\text{Borel}}(\ell, t)$,
it also follows 
 that
$U(\ell)\ \cal{C}_{\text{Borel}}(\ell, |t|\leq z)\subseteq \cal{C}_{\text{Borel}}(\ell, |t|\leq z)$.

\medskip
(vi)  Let 
$M' \in U'(\ell)$ 
be such that
 its  diagonal elements are all equal to some $a\in (\Z/\ell \Z)^{\times}$.
 Let
$M\in {\cal{C}}_{\text{Borel}}(\ell, 0)$. 
 Then $M' M\in B(\ell)$ and $\tr (M' M)=a \tr M=0$, which implies that
 $U'(\ell) \ \cal{C}_{\text{Borel}}(\ell, 0) \subseteq  \cal{C}_{\text{Borel}} (\ell, 0)$. 
\end{proof}

\begin{lemma}\label{gl_2-conj-prop}
Any element of $\GL_2(\ell)$ with a split characteristic polynomial in $\F_{\ell}[X]$ is conjugate to an element of
$\mathcal{B}_{\GL_2}(\ell)$ by  an element of $\SL_2(\ell) := \{M \in \GL_2(\ell): \det M = 1\}$. 
\end{lemma}

\begin{proof}
Let $M \in \GL_2(\ell)$.
Observe that $\lambda_1(M), \lambda_2(M)\in \F_{\ell}^{\times}$. This means that there exists $N \in \GL_2(\ell)$ such that
$$
N M N^{-1}
=
\begin{pmatrix}
\lambda_1(M) & \ast
\\
0 & \lambda_2(M)
\end{pmatrix}\in \mathcal{B}_{\GL_2}(\ell).
$$
If $\det N=1$, we are done. Otherwise, 
by taking 
$$
N' := 
\begin{pmatrix}
(\det N)^{-1} & 0
\\
0 & 1
\end{pmatrix}
N,
$$
we have that $N' M N'^{-1}\in \mathcal{B}_{\GL_2}(\ell)$, as well as that $\det N'=1$, and again we are done.
\end{proof}

\begin{lemma}\label{propr-conj-sets-second}
Every conjugacy class in $\cal{C}(\ell, t)$  or $\cal{C}(\ell, |t|\leq z)$  contains an element of $B(\ell)$.
\end{lemma}
\begin{proof}
Let $M = (M_1, \ldots, M_g) \in {\cal{C}}(\ell, t)$.
Upon fixing $1 \leq i \leq g$, 
we observe that, by Lemma \ref{gl_2-conj-prop}, there exists $N_i\in \SL_2(\ell)$ such that $N_i M_i N_i^{-1}\in \mathcal{B}_{\GL_2}(\ell)$. 
Let  $N:= (N_1, \ldots, N_g)$. Then $N\in G(\ell)$ and 
$N M N^{-1}\in  \mathcal{B}_{\GL_2}(\ell)^g \cap G(\ell) =B(\ell)$. 
We deduce that $M$ is conjugate to an element of $B(\ell)$.
Thus,
every conjugacy class in $\cal{C}(\ell, t)$  contains an element of $B(\ell)$.
 The similar statement about $\cal{C}(\ell, |t|\leq z)$
follows upon noting that every conjugacy class in this set is a conjugacy class in 
$\cal{C}(\ell, t)$ for some  integer $t$ with $|t| \leq z$.
\end{proof}

\begin{lemma}\label{cardin-conj-sets}

\noindent
\begin{enumerate}
\item[(i)]
$|{\cal{C}}_{\text{Torus}}(\ell, t)| \leq 2 (\ell -1)^g$.\\
\item[(ii)]
$|{\cal{C}}_{\text{Borel}}(\ell, t)| = \ell^g |{\cal{C}}_{\text{Torus}}(\ell, t)|\leq 2 (\ell -1)^g\ell^g.$\\
\item[(iii)]
$|\widehat{\cal{C}}_{\text{Borel}}(\ell, t)| =
|{\cal{C}}_{\text{Torus}}(\ell, t)| \leq 2 (\ell -1)^g.$\\
\item[(iv)]
$|\widehat{\cal{C}}'_{\text{Borel}}(\ell, 0)|= 
\ds\frac{|\widehat{\cal{C}}_{\text{Borel}}(\ell, 0)| }{\ell-1}\leq 2 (\ell -1)^{g-1}.$\\
\item[(v)]$|\widehat{\cal{C}}_{\text{Borel}}(\ell, |t|\leq z) |<  5 (\ell-1)^g z$.
\end{enumerate}
\end{lemma}
\begin{proof}
(i) Observe that 
\[
{\cal{C}}_{\text{Torus}}(\ell, t) =
\{(M_1, \ldots, M_g)\in T(\ell): \sum_{1\leq i\leq g} \tr M_i \equiv -t (\mod \ell)\}.
\]
Then
\begin{align*}
|{\cal{C}}_{\text{Torus}}(\ell, t)| 
& =
\ds\sum_{ d \in (\Z/\ell \Z)^{\times}}
\quad
\ds\sum_{
M_1 , \ldots, M_{g-1} \in \mathcal{T}_{\GL_2}(\ell)
\atop{
\det M_1 = \ldots = \det M_g = d
}
}
\# \{
M_g\in \mathcal{T}_{\GL_2}(\ell): \det M_g=d, \tr M_g\equiv -t-\sum_{1\leq i\leq g-1} \tr M_i (\mod \ell)
\}
\\
& \leq \ds\sum_{ d \in (\Z/\ell \Z)^{\times}}
\quad
\ds\sum_{
M_1 , \ldots, M_{g-1} \in \mathcal{T}_{\GL_2}(\ell)
\atop{
\det M_1 = \ldots = \det M_g = d
}
} 2\\
& = 2 (\ell-1)^{g}.
\end{align*}

\medskip
(ii) Observe that 
\[
{\cal{C}}_{\text{Borel}}(\ell, t) =
\{(M_1, \ldots, M_g)\in B(\ell): \sum_{1\leq i\leq g} \tr M_i \equiv -t (\mod \ell)\}.
\]
Then
\begin{align*}
|{\cal{C}}_{\text{Borel}}(\ell, t)| 
& =
\ds\sum_{ d \in (\Z/\ell \Z)^{\times}}
\quad
\ds\sum_{
M_1 , \ldots, M_{g-1} \in \mathcal{T}_{\GL_2}(\ell)
\atop{
\det M_1 = \ldots = \det M_g = d
}
}
\# \{
M_g\in \mathcal{B}_{\GL_2}(\ell): \det M_g=d, \tr M_g\equiv -t-\sum_{1\leq i\leq g-1} \tr M_i (\mod \ell)
\}. 
\end{align*}
Since we have the equality
\begin{align*}
& \# \{
M_g\in \mathcal{B}_{\GL_2}(\ell): \det M_g=d, \tr M_g\equiv -t-\sum_{1\leq i\leq g-1} \tr M_i (\mod \ell)
\} \\
= & 
\ \ell^g \ \# \{
M_g\in \mathcal{T}_{\GL_2}(\ell): \det M_g=d, \tr M_g\equiv -t-\sum_{1\leq i\leq g-1} \tr M_i (\mod \ell)
\} ,
\end{align*}
we obtain that
$|{\cal{C}}_{\text{Borel}}(\ell, t)|  =\ell^g |{\cal{C}}_{\text{Torus}}(\ell, t)|$.

\medskip
(iii) Observe that,
upon identifying $B(\ell)/U(\ell)$ with $T(\ell)$ via Lemma \ref{B/U-lemma},
 there is a bijection between the sets $\widehat{\cal{C}}_{\text{Borel}}(\ell, t)$ and
$\cal{C}_{\text{Torus}}(\ell, t)$. Thus $|\widehat{\cal{C}}_{\text{Borel}}(\ell, t)| =
|{\cal{C}}_{\text{Torus}}(\ell, t)|.$

\medskip
(iv) Since $U(\ell)\leq U'(\ell)$
and
$U'(\ell)/U(\ell) \simeq (\Z/\ell\Z)^{\times}$,
there is a natural surjection
$
B(\ell)/U(\ell)\to B(\ell)/U'(\ell) 
$
whose kernel $U'(\ell)/U(\ell)$ has cardinality $\ell - 1$.
We deduce that $\widehat{\cal{C}}'_{\text{Borel}}(\ell, 0)$ is the image of 
$\widehat{\cal{C}}_{\text{Borel}}(\ell, 0)$ in $B(\ell)/U'(\ell)$. 

By a similar argument as that in the proof of part (vi) of Lemma \ref{propr-conj-sets-first}, 
we obtain that
$U'(\ell)\widehat{\cal{C}}_{\text{Borel}}(\ell, 0)=\widehat{\cal{C}}_{\text{Borel}}(\ell, 0)$.
Thus $\widehat{\cal{C}}_{\text{Borel}}(\ell, 0)$ is
the full inverse image of $\widehat{\cal{C}}'_{\text{Borel}}(\ell, 0)$ in $B(\ell)/U(\ell)$. 

From the above, we infer that
$|\widehat{\cal{C}}'_{\text{Borel}}(\ell, 0)| =
\ds\frac{|\widehat{\cal{C}}_{\text{Borel}}(\ell, 0)| }{\ell-1}\leq 2 (\ell -1)^{g-1}\ell^g$, as claimed.

\medskip
(v) Recall that $$\ds\cal{C}_{\text{Borel}}(\ell, |t|\leq z) =\bigcup_{\substack{t\in \Z\\|t|\leq z}}\cal{C}_{\text{Borel}}(\ell, t).$$
Then, by part (v) of Lemma \ref{propr-conj-sets-first}, we deduce that
\[
\widehat{\cal{C}}_{\text{Borel}}(\ell, |t|\leq z) =\bigcup_{\substack{t\in \Z\\|t|\leq z}}\widehat{\cal{C}}_{\text{Borel}}(\ell, t).
\]
Using part (iii) of the present lemma, we deduce further that
\[
|\widehat{\cal{C}}_{\text{Borel}}(\ell, |t|\leq z)|\leq 2 \ \#\{t\in \Z: |t|\leq z\}(\ell-1)^g < 5z(\ell-1)^g. \]
\end{proof}

\bigskip
\section{Proof of Theorem \ref{main-thm1}}

\subsection{Setting and basic properties}

Let $A$ be an abelian variety defined over $\Q$
that is isogenous over $\Q$ to 
a product of $g$  elliptic curves defined over $\Q$, 
pairwise 
non-isogenous over $\overline{\Q}$
and
each without complex multiplication.
In what follows, we record a few properties of $A$ that we will need in the proofs of our main results.
We keep all  the associated notation introduced in the prior sections, some of which we recall below.
For background on abelian varieties, we refer the reader to 
\cite{Fa83}, \cite{Ho68}, \cite{La83}, \cite{Oo08}, \cite{SeTa68}, and \cite{Wat69}.

We fix a $g$-tuple of elliptic curves, $E_1/\Q$, $\ldots$, $E_g/\Q$,
equipped with a $\Q$-isogeny between $A$ and $E_1 \times \ldots \times E_g$
of minimal degree, denoted $d_A$, among all such choices.
We denote by $N_{E_1}, \ldots, N_{E_g}$ the conductors of $E_1, \ldots, E_g$, respectively,
 and remark that the conductor $N_A$ of $A$ and the product $N_{E_1} \ldots N_{E_g}$ have the same prime factors.


We denote by
$$\rho_A : \Gal\left(\overline{\Q}/\Q\right) \longrightarrow \GL_{2g}(\widehat{\Z})$$
the adelic Galois representation  of $A$, defined by the inverse limit of the residual representations
$$\bar{\rho}_{A, m} : \Gal(\overline{\Q}/\Q) \longrightarrow \GL_{2g}(\Z/m \Z)$$
of the absolute Galois group $\Gal(\overline{\Q}/\Q)$ acting on the $m$-torsion group $A[m] \subset A(\overline{\Q})$, with $m$ a positive integer.
For any prime $\ell$, we denote by 
$$\rho_{A, \ell} : \Gal\left(\overline{\Q}/\Q\right) \longrightarrow \GL_{2g}(\Z_{\ell})$$
the $\ell$-adic representation of $\Gal(\overline{\Q}/\Q)$ on the  $\ell$-adic Tate module ${\displaystyle \lim_{\leftarrow \atop{n}} A[\ell^n}]$. 

Similarly to  notation (\ref{group-G-of-ell}) introduced in Section 4, 
we set
$$
G(\widehat{\Z})
:= \{(M_1, \ldots, M_g)\in  \GL_2(\widehat{\Z})^g: \det M_1 = \ldots = \det M_g\},
$$
 $$
G(m)
:= \{(M_1, \ldots, M_g)\in  \GL_2(\Z/m \Z)^g: \det M_1 = \ldots = \det M_g\},
$$
and recall that 
\begin{equation*}\label{image-galois-av-absolute}
\Im \rho_A \subseteq G(\widehat{\Z}),
\end{equation*}
\begin{equation*}\label{image-galois-av-residual}
\Im \bar{\rho}_{A, m} \subseteq G(m),
\end{equation*}
where $m$ is an arbitrary  positive integer.

For any prime $p \nmid N_A$, we denote by $\Frob_p$ the Frobenius at $p$ in $\Gal(\overline{\Q}/\Q)$,
and  for any prime $\ell \neq p$, we denote by
$$P_{A, p}(X) := \det (X I_{2g} - \rho_{A, \ell}(\Frob_p))$$
the characteristic polynomial of $\rho_{A, \ell}(\Frob_p)$ in $\GL_{2g}(\Z_{\ell})$.

We recall that
$P_{A, p}(X)$ is the characteristic polynomial of the Frobenius endomorphism 
 acting on the reduction of $A$ modulo $p$. 
We also recall
 that $P_{A, p}(X)$ is a $p$-Weil polynomial,
has integer coefficients  independent of $\ell$, 
and satisfies the congruence
\begin{equation}\label{p-Weil-char-m}
P_{A, p}(X) \equiv \det \left(X I_{2g} - \bar{\rho}_{A, m} \left(\left(\frac{\Q(A[m])/\Q}{p}\right)\right)  \right) (\mod m)
\end{equation}
for any integer $m$ coprime to  $p$. 

In analogy with the discussion and notation for $A$,
for each $1 \leq i \leq  g$
we denote by
$\rho_{E_i}$, $\overline{\rho}_{E_i, m}$, and $\rho_{E_i, \ell}$ the 
adelic, residual, and $\ell$-adic Galois representations of $E_i$,
and
for each prime $p \nmid N_{E_i}$
we denote  by
$P_{E_i, p}(X)$
the associated $p$-Weil polynomial.
We recall 
the congruence relation
\begin{equation}\label{p-Weil-char-m-ec}
P_{E_i, p}(X) \equiv 
\det \left(X I_{2} - \bar{\rho}_{E_i, m} \left(\frac{\Q(E_i[m])/\Q}{p}\right)  \right)(\mod m),
\end{equation}
which holds for any integer $m$ coprime to  $p$. 

Now let us fix a prime $p \nmid N_A d_A$.
We write the polynomials $P_{A, p}(X)$, $P_{E_i, p}(X)$ explicitly as 
\begin{eqnarray*}
P_{A, p}(X) 
&=& 
 X^{2g} + a_{1, p}(A) X^{2g-1} + \ldots  + a_{2g-1, p}(A) X + a_{2g, p}(A)  \in \Z[X],
\\
P_{E_i, p}(X)
&=& X^2 - a_p(E_i) X + p  \in \Z[X].
\end{eqnarray*}
On one hand, 
 we have the  polynomial relation
\begin{equation}\label{polyn-product}
P_{A, p}(X) = P_{E_1, p}(X) \cdot \ldots \cdot P_{E_g, p}(X),
\end{equation}
from which we obtain that
\begin{equation}\label{Frob-trace-sum}
a_{1, p}(A) = - (a_p(E_1) + \ldots + a_p(E_g)).
\end{equation}
On the other hand,  by restricting the residual representations $\overline{\rho}_{A, m}$, $\overline{\rho}_{E_i, m}$
to the division fields $\Q(A[m])$, $\Q(E_i[m])$, respectively,
we deduce from (\ref{p-Weil-char-m}) - (\ref{p-Weil-char-m-ec}) that,
 for any integer $m$ coprime to  $p$,
we have
\begin{eqnarray}
\tr  \bar{\rho}_{A, m} \left( \left(\frac{\Q(A[m])/\Q}{p}\right)  \right)  
&\equiv&
- a_{1, p}(A) (\mod m),
\label{trace-congruence-av}
\\
\tr  \bar{\rho}_{E_i, m} \left( \left(\frac{\Q(E_i[m])/\Q}{p}\right)  \right)  
&\equiv&
a_{p}(E_i) (\mod m).
\label{trace-congruence-ec}
\end{eqnarray}

\noindent
From (\ref{Frob-trace-sum}), (\ref{trace-congruence-av}) and (\ref{trace-congruence-ec}), we deduce that,
for any integer $m$ coprime to $p$, we have
\begin{equation}\label{trace-sum-congruence}
\tr  \bar{\rho}_{A, m} \left( \left(\frac{\Q(A[m])/\Q}{p}\right)  \right)  
\equiv
-
\ds\sum_{1 \leq i \leq g}
\tr  \bar{\rho}_{E_i, m} \left( \left(\frac{\Q(E_i[m])/\Q}{p}\right)  \right) 
(\mod m).
\end{equation}


\bigskip
\subsection{A conditional reduction step}

For an abelian variety $A/\Q$, a prime $p \nmid N_A$, and an integer $t$,
we remark that,
congruence (\ref{trace-congruence-av}) already allows us 
to relax the equality $a_{1, p}(A) = t$ to
the congruence $\tr  \bar{\rho}_{A, \ell} \left( \left(\frac{\Q(A[\ell])/\Q}{p}\right)  \right)   \equiv  - t (\mod \ell)$  for some prime $\ell$
and then to apply an effective version of the Chebotarev Density Theorem
 in the extension $\Q(A[\ell])/\Q$.
 However,
we can obtain better results by relaxing the equality $a_{1, p}(A) = t$  to a Chebotarev condition
that takes place 
in the extension $\Q(A[\ell])^{U(\ell)}/\Q(A[\ell])^{B(\ell)}$. 
The key ingredient is as follows.

\begin{lemma}\label{lemma-max-pi-A-x-ell-t}
Let $t \in \Z$.
Let $A/\Q$ be an abelian variety defined over $\Q$
that is isogenous over $\Q$ to a product $E_1 \times \ldots \times E_g$
of elliptic curves $E_1/\Q$, $\ldots$, $E_g/\Q$, without complex multiplication
and pairwise non-isogenous over $\overline{\Q}$.
 For each $x>2$, let $y = y(x) > 2$, $u =u(x) > 2$ be such that
\begin{equation}\label{hypothesis3-applic}
u \leq y.
\end{equation}
Assume that, for any $\varepsilon > 0$,
\begin{equation}\label{hypothesis4-applic}
u \geq c'(\varepsilon)
y^{\frac{1}{2}} (\log y)^{2 + \varepsilon} 
\
\text{for some} \ c'(\varepsilon) > 0
\end{equation}
and
\begin{equation}\label{hypothesis5-applic}
\ds\lim_{x \rightarrow \infty} \frac{\log x}{(\log y)^{1 + \varepsilon}} = 0.
\end{equation}
Assume GRH for Dedekind zeta functions. 
 Then,  
  for any $\varepsilon > 0$,
 there exists a constant $c(\varepsilon) > 0$ such that,
 for any positive real number $z$
 and for any sufficiently large $x$, we have
\begin{equation}\label{max-lemma-bound-applic}
\pi_A(x, t)
\leq
c(\varepsilon)
\ds\max_{y \leq \ell \leq y + u} 
\pi_A(x, \ell, t)
\end{equation}
and
\begin{equation}\label{max-lemma-bound-applic-coro}
\sum_{\substack{t\in \Z\\|t|\leq z}}\ \pi_A(x, t)
\leq
c(\varepsilon)
\ds\max_{y \leq \ell \leq y + u} 
\sum_{\substack{t\in \Z\\|t|\leq z}}\ \pi_A(x, \ell, t),
\end{equation}
where, as before,
\begin{equation*}
\pi_A(x, t) = \#\{p \leq x: p \nmid N_A, a_{1, p}(A) = t\},
\end{equation*}
and
where
\begin{eqnarray}\label{pi-A-x-ell-t-def}
\pi_A(x, \ell, t)
&:=&
\#\left\{p \leq x: p \nmid \ell N_A, a_{1, p}(A) = t,
\right.
\nonumber
\\
&&
\hspace*{1.6cm}
\left.
\ell \ \text{splits completely in each of the fields} \ \Q(\pi_p(E_1)), \ldots,  \Q(\pi_p(E_g))
\right\}.
\nonumber
\end{eqnarray}
\end{lemma}
\begin{proof}
We apply Lemma \ref{max-lemma} to
the set ${\cal{S}} = \left\{p:  p \nmid N_A, a_{1, p}(A) = t \right\}$,
the field $K_p = $   the compositum of the fields $\Q(\pi_p(E_1)), \ldots, \Q(\pi_p(E_g))$,
the conjugacy class ${\cal{C}}_p = \left\{\id_{K_p}\right\}$,
and the function $f(v) = \log v$.

Note that $K_p/\Q$ is Galois and that $1 \neq [K_p:\Q] \mid 2^g$.
Also note that if a prime $\ell$ is ramified in $K_p$, then
$\ell \mid \left(4p - a_p(E_1)^2\right) \ldots \left(4p - a_p(E_g)^2\right)$ (see \cite[Theorem, p. 43]{To55}).
Therefore, 
upon recalling  (\ref{hensel}),
we have that
$\ds \log |d_{K_p}| \ll_{E_1, \ldots, E_g} \log p$.
 Since all hypotheses of Lemma \ref{max-lemma} are satisfied, bound (\ref{max-lemma-bound-applic}) follows. 
 Bound (\ref{max-lemma-bound-applic-coro}) follows from a similar argument, except for 
 taking
  ${\cal{S}} = \left\{p:  p \nmid N_A, |a_{1, p}(A)| \leq z\right\}$,
  while keeping 
  $K_p$, ${\cal{C}}_p$, and $f$ as before.

\end{proof}

\bigskip
\subsection{Proof of Theorem \ref{main-thm1} for arbitrary $t$}
Let $A$ be an abelian variety defined over $\Q$
that is isogenous over $\Q$ to 
a product of $g$  elliptic curves defined over $\Q$, 
pairwise 
non-isogenous over $\overline{\Q}$
and
each without complex multiplication.
As in Subsection 6.1, 
we fix a $g$-tuple of elliptic curves, $E_1/\Q$, $\ldots$, $E_g/\Q$,
equipped with a $\Q$-isogeny between $A$ and $E_1 \times \ldots \times E_g$
of minimal degree, $d_A$, among all such choices.
We keep all the associated notation introduced so far.

Let $t \in \Z$ and  $x > 2$.
By Lemma \ref{lemma-max-pi-A-x-ell-t},
we can bound $\pi_A(x, t)$ from above if we succeed in bounding  $\pi_A(x, \ell, t)$ from above for some suitably chosen prime $\ell = \ell(x)$.

For now, we fix a prime $\ell \nmid 2 g$ 
such that 
\begin{equation}\label{image-large-ell}
\Im \overline{\rho}_{A, \ell} = G(\ell).
\end{equation}
 The existence of infinitely many such primes $\ell$ is ensured by \cite[Th\'{e}or\`{e}me 6, p. 324]{Se72} and \cite[Theorem 1.1, p. 387]{Lo16} under our hypotheses that   
$E_1/\Q, \ldots, E_g/\Q$
are without complex multiplication
and
pairwise non-isogenous over $\overline{\Q}$.
Indeed, from  \cite[Theorem 1.1, p. 387]{Lo16}, we infer that
 there exists a least positive integer $m_A$  such that  $\Im \rho_A$ 
 equals the inverse image under the canonical projection 
 $G (\widehat{\Z}) \longrightarrow G(m_A)$ of $\Im \overline{\rho}_{A, m_A}$.
 In particular, this means that
$\Im \bar{\rho}_{A, \ell} = G(\ell)$
for any prime $\ell \nmid m_A$.
Later on,  we will choose $\ell$ in an explicit interval that depends on $x$.

The choice of a prime $\ell$ such that (\ref{image-large-ell}) holds 
 allows us to consider the subextensions of $\Q(A[\ell])$ fixed by the groups $U(\ell)$ and $B(\ell)$,
 namely
$$
\Q \subset \Q(A[\ell])^{B(\ell)} \subset \Q(A[\ell])^{U(\ell)} \subset \Q(A[\ell]),
$$
for which we immediately obtain the following Galois group structures:
\begin{eqnarray*}
\Gal(\Q(A[\ell])/\Q) &\simeq& G(\ell),
\\
\Gal\left(\Q(A[\ell])/\Q(A[\ell])^{B(\ell)}\right) &\simeq& B(\ell),
\\
\Gal\left(\Q(A[\ell])/\Q(A[\ell])^{U(\ell)}\right) &\simeq& U(\ell),
\\
\Gal\left(\Q(A[\ell])^{U(\ell)}/\Q(A[\ell])^{B(\ell)}\right) &\simeq& B(\ell)/U(\ell).
\end{eqnarray*}
In relation to these Galois groups, we recall  from part (v) of Lemma \ref{centr-G-B-U}
that $U(\ell) \unlhd B(\ell)$ and that $B(\ell)/U(\ell)$ is abelian.

Using 
(\ref{polyn-product}),
(\ref{trace-sum-congruence}),
and
(\ref{image-large-ell}),
we deduce that
every prime 
$p \nmid d_A$
counted by $\pi_A(x, \ell, t)$ 
is also  counted  by $\pi_{\cal{C}(\ell, t)}(x, \Q(A[\ell])/\Q)$,
where $\cal{C}(\ell, t)$ is 
the union of conjugacy classes of $G(\ell)$ 
studied in Lemma \ref{propr-conj-sets-first} of Section 5. 
Thus
\begin{equation}\label{pi-ell-t-C-ell-t}
\pi_A(x, \ell, t) \ll_A \pi_{\cal{C}(\ell, t)}(x, \Q(A[\ell])/\Q).
\end{equation}

To estimate $\pi_{\cal{C}(\ell, t)}(x, \Q(A[\ell])/\Q)$, we invoke part (ii) of 
Corollary \ref{cheb-upper-bound-as-wanted} of Section 2 with
$K = \Q$, $L = \Q(A[\ell])$, $G = G(\ell)$, $H = B(\ell)$, $N = U(\ell)$, and ${\cal{C}} = {\cal{C}}(\ell, t)$.
The hypotheses of this corollary hold thanks to our assumption of GRH 
and 
to part (v) of Lemma \ref{centr-G-B-U},
part (v) of Lemma \ref{propr-conj-sets-first},
and Lemma \ref{propr-conj-sets-second}.
We deduce that
\begin{eqnarray*}
\pi_{\cal{C}(\ell, t)}(x, \Q(A[\ell])/\Q)
&\ll&
\frac{\left|\widehat{\cal{C}}_{\text{Borel}}(\ell, t)\right| \cdot |U(\ell)|}{|B(\ell)|} \cdot \frac{x}{\log x}
\\
&+&
\left|\widehat{{\cal{C}}}_{\text{Borel}}(\ell, t)\right|^{\frac{1}{2}} [\Q(A[\ell])^{B(\ell)} :\Q] \frac{x^{\frac{1}{2}}}{\log x} \log M\left(\Q(A[\ell])^{U(\ell)}/\Q(A[\ell])^{B(\ell)}\right)
\\
&+&
\left(\frac{x^{\frac{1}{2}}}{\log x} + \log M \left( \Q(A[\ell])/\Q\right) \right)
\\
&+&
[\Q(A[\ell])^{B(\ell)}:\Q] \left(\frac{x^{\frac{1}{2}}}{\log x} 
+ \log M\left(\Q(A[\ell])^{U(\ell)}/\Q(A[\ell])^{B(\ell)}\right)\right).
\end{eqnarray*}
Using Lemma  \ref{counting-G(l)-etc}  and Lemma \ref{cardin-conj-sets}, we deduce from the above that
\begin{equation}\label{interm-bound}
\pi_{\cal{C}(\ell, t)}(x, \Q(A[\ell])/\Q)
\ll
\frac{x}{\ell \log x}
+
 (\ell-1)^{\frac{g}{2}} (\ell + 1)^g
\frac{x^{\frac{1}{2}}}{\log x} \log M\left(\Q(A[\ell])^{U(\ell)}/\Q(A[\ell])^{B(\ell)}\right).
\end{equation}

We estimate 
$ \log M\left(\Q(A[\ell])^{U(\ell)}/\Q(A[\ell])^{B(\ell)}\right)$
by putting together (\ref{hensel}),
Lemma \ref{counting-G(l)-etc},
and 
 the N\'{e}ron-Ogg-Shafarevich criterion for $A/\Q$, which states that the extension $\Q(A[\ell])/\Q$ is unramified outside of 
$\ell N_A$  (see \cite[Theorem 1, p. 493]{SeTa68}).
We obtain that
\begin{equation}\label{long-log-M-t}
\log M\left(\Q(A[\ell])^{U(\ell)}/\Q(A[\ell])^{B(\ell)}\right)
\leq
2 \log |B(\ell)/U(\ell)| + 2\log (\ell N_A) + \log 2
\ll
g \log (\ell N_A).
\end{equation}

Using this bound in the upper-estimate (\ref{interm-bound}) for $\pi_{\cal{C}(\ell, t)}(x, \Q(A[\ell])/\Q)$, we deduce that
\begin{equation}\label{remark-t}
\pi_{\cal{C}(\ell, t)}(x, \Q(A[\ell])/\Q)
\ll
\frac{x}{\ell \log x}
+
g\  \ell^{\frac{3g}{2}}
\frac{x^{\frac{1}{2}}}{\log x} \log (\ell N_A).
\end{equation}

Finally, we apply Lemma \ref{lemma-max-pi-A-x-ell-t} with
$$y \asymp \frac{x^{\frac{1}{3 g + 2}}}{(\log x)^{\frac{2}{3 g + 2}}}$$
and $$u \asymp y^{\frac{1}{2}} (\log y)^{2 + \varepsilon}$$ for an arbitrarily chosen $\varepsilon > 0$,
and conclude that
$$
\pi_A(x, t)
\leq
c(\varepsilon)
\ds\max_{y \leq \ell \leq y + u} 
\pi_A(x, \ell, t)
\ll_{\varepsilon, A}
\frac{
x^{1 - \frac{1}{3g+2} }
}{
(\log x)^{1 - \frac{2}{3 g + 2} }}.
$$

\bigskip
\subsection{Proof of Theorem \ref{main-thm1} for  $t = 0$}
We keep the setting of Subsection 6.3.
 Let $x > 2$.
As before,  we take $\ell \nmid 2 g$ to be a sufficiently large prime such that 
(\ref{image-large-ell})
holds.
Then taking $t = 0$ in (\ref{pi-ell-t-C-ell-t}), 
we have that
\begin{equation}\label{pi-ell-t-C-ell-0}
\pi_A(x, \ell, 0) \ll_A \pi_{\cal{C}(\ell, 0)}(x, \Q(A[\ell])/\Q).
\end{equation}
To estimate  $\pi_{\cal{C}(\ell, 0)}(x, \Q(A[\ell])/\Q)$,
we follow the strategy of the previous subsection, the main difference being that we will now work
with the subextensions of $\Q(A[\ell])$ fixed by the groups $U'(\ell)$ and $B(\ell)$, 
instead of by $U(\ell)$ and $B(\ell)$.
More precisely, we will work with the extensions
$$
\Q \subset \Q(A[\ell])^{B(\ell)} \subset \Q(A[\ell])^{U'(\ell)} \subset \Q(A[\ell]).
$$

Observe that 
$\Gal\left(\Q(A[\ell])/\Q(A[\ell])^{U'(\ell)}\right) \simeq U'(\ell)$
and
$\Gal\left(\Q(A[\ell])^{U'(\ell)}/\Q(A[\ell])^{B(\ell)}\right) \simeq B(\ell)/U'(\ell)$,
and
recall from part (vi) of Lemma \ref{centr-G-B-U} that 
$U'(\ell) \unlhd B(\ell)$ and that $B(\ell)/U'(\ell)$ is abelian. 
Upon also recalling part (vi) of Lemma \ref{propr-conj-sets-first}
and Lemma \ref{propr-conj-sets-second},
we apply 
part (ii) of 
Corollary \ref{cheb-upper-bound-as-wanted} of Section 2 with
$K = \Q$, $L = \Q(A[\ell])$, $G = G(\ell)$, $H = B(\ell)$, $N = U'(\ell)$, and ${\cal{C}} = {\cal{C}}(\ell, 0)$,
and deduce that
\begin{eqnarray*}
\pi_{\cal{C}(\ell, 0)}(x, \Q(A[\ell])/\Q)
&\ll&
\frac{\left|\widehat{\cal{C}}'_{\text{Borel}}(\ell, 0)\right| \cdot |U'(\ell)|}{|B(\ell)|} \cdot \frac{x}{\log x}
\\
&+&
\left|\widehat{{\cal{C}}}'_{\text{Borel}}(\ell, 0)\right|^{\frac{1}{2}} [\Q(A[\ell])^{B(\ell)} :\Q] \frac{x^{\frac{1}{2}}}{\log x} \log M\left(\Q(A[\ell])^{U'(\ell)}/\Q(A[\ell])^{B(\ell)}\right)
\\
&+&
\left(\frac{x^{\frac{1}{2}}}{\log x} + \log M\left(\Q(A[\ell])/\Q\right)\right)
\\
&+&
[\Q(A[\ell])^{B(\ell)}:\Q] \left(\frac{x^{\frac{1}{2}}}{\log x} 
+ \log M\left(\Q(A[\ell])^{U'(\ell)}/\Q(A[\ell])^{B(\ell)}\right)\right).
\end{eqnarray*}
Using Lemma  \ref{counting-G(l)-etc}  and Lemma \ref{cardin-conj-sets}, we infer from the above that
\begin{eqnarray*}
\pi_{\cal{C}(\ell, 0)}(x, \Q(A[\ell])/\Q)
&\ll&
\frac{x}{\ell \log x}
+
 \ell^{\frac{g-1}{2}} (\ell + 1)^g
\frac{x^{\frac{1}{2}}}{\log x} \log M(\Q(A[\ell])^{U'(\ell)}/\Q(A[\ell])^{B(\ell)}).
\end{eqnarray*}
Similarly to how we deduced (\ref{long-log-M-t}),
we obtain that
\begin{equation}\label{long-log-M-0}
\log M\left(\Q(A[\ell])^{U'(\ell)}/\Q(A[\ell])^{B(\ell)}\right)
\leq
2 \log |B(\ell)/U'(\ell)| +2 \log (\ell N_A) + \log 2
\ll
g \log (\ell N_A).
\end{equation}
Thus,
\begin{eqnarray*}
\pi_{\cal{C}(\ell, t)}(x, \Q(A[\ell])/\Q)
&\ll&
\frac{x}{\ell \log x}
+
g\ \ell^{\frac{3g-1}{2}} 
\frac{x^{\frac{1}{2}}}{\log x} \log (\ell N_A).
\end{eqnarray*}

Now we apply Lemma \ref{lemma-max-pi-A-x-ell-t} with
$$y \asymp \frac{x^{\frac{1}{3 g + 1}}}{(\log x)^{\frac{2}{3 g + 1}}}$$
and $$u \asymp y^{\frac{1}{2}} (\log y)^{2 + \varepsilon}$$ for an arbitrarily chosen $\varepsilon > 0$.
We conclude that
$$
\pi_A(x, 0)
\leq
c(\varepsilon)
\ds\max_{y \leq \ell \leq y + u} 
\pi_A(x, \ell, 0)
\ll_{\varepsilon, A}
\frac{
x^{1 - \frac{1}{3g+1} }
}{
(\log x)^{1 - \frac{2}{3 g + 1} }}.
$$
This completes the proof of Theorem \ref{main-thm1}.

\bigskip
\subsection{Proof of applications (\ref{corollary-1}) and (\ref{corollary-2})}

We keep the setting of Subsection 6.3.
Let $t \in \Z$ and  $x > 2$.
By applying Theorem \ref{main-thm1}, we deduce right away that
\begin{eqnarray*}
\#\left\{p \leq x: p \nmid N_A, a_{1, p}(A) \neq t \right\}
=
\pi(x)
-
\#\left\{p \leq x: p \nmid N_A, a_{1, p}(A) = t\right\}
-
\#\left\{p \leq x: p \mid N_A\right\}
\sim
\pi(x),
\end{eqnarray*}
which confirms (\ref{corollary-1}).
To confirm 
 (\ref{corollary-2}),
 we fix $\varepsilon > 0$
 and
 make the following  observations, 
 the third of which uses Theorem \ref{main-thm1}:
\begin{eqnarray*}
\pi(x)
&=&
\#\left\{p \leq x: p \nmid N_A, |a_{1, p}(A)| > p^{\frac{1}{3g+2} - \varepsilon}  \right\}
+
\#\left\{p \leq x: p \nmid N_A, |a_{1, p}(A)| \leq p^{\frac{1}{3g+2} - \varepsilon}  \right\}
+
\#\left\{p \leq x: p \mid N_A\right\}
\\
&=&
\#\left\{p \leq x: p \nmid N_A, |a_{1, p}(A)| > p^{\frac{1}{3g+2} - \varepsilon} \right\}
+
\O\left(
\ds\sum_{t \in \Z \atop{|t| \leq x^{\frac{1}{3g+2} - \varepsilon} }}
\pi_A(x, t)
\right)
+
\O_{A}(1)
\\
&=&
\#\left\{p \leq x: p \nmid N_A, |a_{1, p}(A)| > p^{\frac{1}{3g+2} - \varepsilon} \right\}
+
\O_{A, \varepsilon}\left(
x^{\frac{1}{3g+2} - \varepsilon}
\cdot
\ds\frac{x^{1 - \frac{1}{3 g + 2}}}{ (\log x)^{1 - \frac{2}{3 g + 2}} }
\right)
\\
&=&
\#\left\{p \leq x: p \nmid N_A, |a_{1, p}(A)| > p^{\frac{1}{3g+2} - \varepsilon} \right\}
+
\O_{A, \varepsilon}\left(
\ds\frac{x^{1 - \varepsilon} }{ (\log x)^{1 - \frac{2}{3 g + 2}} }
\right)
\\
&=&
\#\left\{p \leq x: p \nmid N_A, |a_{1, p}(A)| > p^{\frac{1}{3g+2} - \varepsilon} \right\}
+
\o\left(\pi(x)\right).
\end{eqnarray*}

\bigskip
\subsection{Proof of Theorem \ref{main-thm2}}
We keep the setting of Subsection 6.3.
Let $x > 2$ and $z = z(x) > 1$.  
Note that, for any prime $\ell$, we have
\begin{align}\label{trace-sum}
\sum_{\substack{t\in \Z \\|t|\leq z}}\
\pi_A(x, \ell, t)
&:=
\#\left\{p \leq x: p \nmid \ell N_A, |a_{1, p}(A)| \leq z,
\right.
\nonumber
\\
&
\hspace*{1.6cm}
\left.
\ell \ \text{splits completely in each of the fields} \ \Q(\pi_p(E_1)), \ldots,  \Q(\pi_p(E_g))
\right\}.
\end{align}


As before,  
observe that 
every prime 
$p \nmid d_A$
 which is counted on
the right-hand side of 
 (\ref{trace-sum}) 
 is also counted
  by $\pi_{\cal{C}(\ell, |t|\leq z) }(x, \Q(A[\ell])/\Q)$.
Thus
\[
\ds\sum_{\substack{t\in \Z\\|t|\leq z}} \pi_A(x, \ell, t) \ll_A \pi_{\cal{C}(\ell, |t|\leq z) }(x, \Q(A[\ell])/\Q).
\]

We take $\ell \nmid 2 g$ to be a sufficiently large prime such that 
(\ref{image-large-ell})
holds.
Recalling part (vi) of Lemma \ref{propr-conj-sets-first}
and Lemma \ref{propr-conj-sets-second},
we apply 
part (ii) of 
Corollary \ref{cheb-upper-bound-as-wanted} of Section 2 with
$K = \Q$, $L = \Q(A[\ell])$, $G = G(\ell)$, $H = B(\ell)$, $N = U(\ell)$, and ${\cal{C}} = {\cal{C}}(\ell, |t|\leq z)$,
and  obtain
\begin{eqnarray*}\label{bound-cor}
\pi_{\cal{C}(\ell, |t|\leq z) }(x, \Q(A[\ell])/\Q)
&\ll&
\frac{\left|\widehat{\cal{C}}_{\text{Borel}}(\ell, |t|\leq z) \right|
}{\left|B(\ell)/U(\ell)\right|
} \cdot \frac{x}{\log x}
\\
&+&
\left|\widehat{\cal{C}}_{\text{Borel}}(\ell, |t|\leq z) \right|^{\frac{1}{2}}
\left[\Q(A[\ell])^{B(\ell)} : \Q\right] \frac{x^{\frac{1}{2}}}{\log x}
\left(\log M\left(\Q(A[\ell])^{U(\ell)}/\Q(A[\ell])^{B(\ell)}\right)\right).
\nonumber
\end{eqnarray*}
Using Proposition  \ref{counting-G(l)-etc}  and Lemma \ref{cardin-conj-sets}, we infer from the above that
\[
\pi_{\cal{C}(\ell, |t|\leq z) }(x, \Q(A[\ell])/\Q)\ll
\frac{z}{\ell}\cdot \frac{x}{\log x} + g \ z^{\frac{1}{2}}\ell^{\frac{3g}{2}} \frac{x^{\frac{1}{2}}}{\log x}\log (\ell N_A).
\]

Applying (\ref{max-lemma-bound-applic-coro}) in Lemma \ref{lemma-max-pi-A-x-ell-t} with 
$y \asymp z(\log x)^{\eta}$ for some fixed arbitrary $\eta>0$,  we obtain that
\begin{align*}
\sum_{\substack{t\in \Z\\|t|\leq z}}\pi_A(x, t)& \ll_{A} \frac{x}{(\log x)^{1+\eta}}+x^{\frac{1}{2}}  (\log x)^{\frac{3g\eta}{2}-1}z^{\frac{3g+1}{2}}(\log z+\log\log x). 
\end{align*}


Now fix an arbitrary $\varepsilon>0$ and set
$
z := x^{\frac{1}{3g+1}-\varepsilon}.
$ 
Then
\[
\sum_{\substack{t\in \Z\\ |t|\leq x^{\frac{1}{3g+1}-\varepsilon}}}\pi_A(x, t)
\ll_{A, \varepsilon} 
\frac{x}{(\log x)^{1+\eta}} = \o(\pi(x)).
\]
Theorem \ref{main-thm2} follows from the above upper bound and the following observations:
\begin{eqnarray*}\pi(x)&=&\#\left\{p \leq x: p \nmid N_A, |a_{1, p}(A)| > p^{\frac{1}{3g+1} - \varepsilon}  \right\}+\#\left\{p \leq x: p \nmid N_A, |a_{1, p}(A)| \leq p^{\frac{1}{3g+1} - \varepsilon}  \right\}+\#\left\{p \leq x: p \mid N_A\right\}\\
&=&\#\left\{p \leq x: p \nmid N_A, |a_{1, p}(A)| > p^{\frac{1}{3g+1} - \varepsilon} \right\}+\ds\sum_{t \in \Z \atop{|t| \leq x^{\frac{1}{3g+1} - \varepsilon} }}\pi_A(x, t)+\O_A(1)\\
&=&\#\left\{p \leq x: p \nmid N_A, |a_{1, p}(A)| > p^{\frac{1}{3g+1} - \varepsilon} \right\}+\o\left(\pi(x)\right).\end{eqnarray*}

\bigskip
\section{Final remarks}

It is a natural question to ask whether one can find 
subgroups and conjugacy classes other than those chosen in Sections 4 - 5
 and use them, following a strategy similar to the one of Theorem \ref{main-thm1},
to obtain a better upper bound for $\pi_A(x, t)$.
We relegate this question to future research.
Below are candidates of subgroups and conjugacy classes which, used mutatis mutandis in our proof,
 do not  improve on Theorem \ref{main-thm1}.

We keep the setting of Subsection 6.3.
We could replace 
$B(\ell)$ with  $T(\ell)$, 
$U(\ell)$ with $\{I_{2g}\}$,  
and  
$\widehat{{\cal{C}}}_{\text{Borel}}(\ell, t)$
with
$\cal{C}_{\text{Torus}}(\ell, t)$.
%
For every $t\in \Z$, this would lead to the upper bound 
\begin{equation}\label{tours-bound}
\pi_A(x, t)\ll_A \frac{x^{1-\frac{1}{5g+2}}}{(\log x)^{1-\frac{2}{5g+2}}}.
\end{equation}

Alternatively, 
we could replace 
$B(\ell)$ with  
$C^{ns}_{\xi}(\ell)^g\cap G(\ell)$,
$U(\ell)$ with $\{I_{2g}\}$,  
and  
$\widehat{{\cal{C}}}_{\text{Borel}}(\ell, t)$
with
$\cal{C}_{ns}(\ell, \xi, t)$,
where 
$C^{ns}_{\xi}(\ell)$ is the non-split Cartan subgroup of $\GL_2(\ell)$ defined by
$$
C^{ns}_{\xi}(\ell) := 
\left\{
\begin{pmatrix}
a & \xi b\\
b & a
\end{pmatrix}: 
a, b\in \F_{\ell}, (a, b)\neq (0, 0)
\right\} 
$$
for some   non-square element $\xi$ of $\F_{\ell}^{\times}$,
and where
$$
\hspace*{-0.5cm}
\cal{C}_{ns}(\ell, \xi, t) :=\left\{
(M_1, \ldots, M_g) \in C^{ns}_{\xi}(\ell)^g\cap G(\ell) :
\ds\sum_{
1 \leq i \leq g
}
\tr M_i
\equiv
- t (\mod \ell)
\right\}.
$$
Note that  $C^{ns}_{\xi}(\ell)^g\simeq( \F_{\ell^2}^{\times})^g$, 
where $\F_{\ell^2}$ denotes the finite field with $\ell^2$ elements.
Hence $C^{ns}_{\xi}(\ell)^g$ is an abelian group.
With the above choices,
instead of considering the counting function $\pi_{A}(x, \ell, t)$, 
we would consider 
\[
\pi^{ns}_{A}(x, \ell, t) 
:=
\#\{p\leq x: p\nmid N_A, a_{1, p}(A)=t, \ \ell \ \text{inert in each of the fields} \ \Q(\pi_p(E_1)), \ldots,  \Q(\pi_p(E_g))\}
\]
and prove 
an upper bound similar to that
of
 Lemma \ref{lemma-max-pi-A-x-ell-t}.
Then,  applying
Corollary \ref{cheb-upper-bound-as-wanted},
we would obtain the previous bound (\ref{tours-bound}), which is weaker than that of Theorem \ref{main-thm1}.

\bigskip

{\small{

}

\end{document}